\documentclass[a4 paper, 12pt]{article} 
\usepackage{amssymb,amsthm}
\usepackage{amsmath}
\usepackage{subcaption}

\usepackage{tikz}
\usepackage{graphicx}
\usepackage[shortlabels]{enumitem}

\graphicspath{{microlocalfigures1/}}
\usetikzlibrary{decorations.pathmorphing}
\tikzset{snake it/.style={decorate, decoration=snake}}

\title{Microlocal analysis of a spindle transform}
\author{James Webber\thanks{Corresponding author, supported by Engineering and Physical Sciences Research Council and Rapiscan systems, CASE studentship.} \and Sean Holman\thanks{Second author supported by Engineering and Physical Sciences Research Council (EP/M016773/1).}}
\date{}
\theoremstyle{plain}
\newtheorem{theorem}{Theorem}

\theoremstyle{definition}
\newtheorem{definition}{Definition}
\newtheorem{proposition}{Proposition}
\newtheorem{corollary}{Corollary}

\begin{document}
\maketitle

\begin{abstract}
An analysis of the stability of the spindle transform, introduced in \cite{WL}, is presented. We do this via a microlocal approach and show that the normal operator for the spindle transform is a type of paired Lagrangian operator with ``blowdown--blowdown" singularities analogous to that of a limited data synthetic aperture radar (SAR) problem studied by Felea et. al. \cite{felea}. We find that the normal operator for the spindle transform belongs to a class of distibutions $I^{p,l}(\Delta\cup\widetilde{\Delta},\Lambda)$ studied by Felea and Marhuenda in \cite{felea,theoremref}, where $\widetilde{\Delta}$ is reflection through the origin, and $\Lambda$ is associated to a rotation artefact. Later, we derive a filter to reduce the strength of the image artefact and show that it is of convolution type. We also provide simulated reconstructions to show the artefacts produced by $\Lambda$ and show how the filter we derived can be applied to reduce the strength of the artefact. 
\end{abstract}

\section{Introduction}
\begin{figure}[!h]
\centering
\begin{tikzpicture}[scale=4.5]
\path ({-cos(30)},0.5) coordinate (S);
\path ({cos(30)}, 0.5) coordinate (D);
\path [rotate=45]  (-0.5,0.26) coordinate (w);
\path (-0.13,1.22) coordinate (a);
\draw [rotate=45] (0,-0.6) circle [radius=1];
\draw  [rotate=45]  (0,0.6) circle [radius=1];
\draw [rotate=45]  [->] (0,0)--(0.35,0);
\draw [rotate=45]  [<->] (0,-0.6)--(0,0);
\draw [rotate=45] (0,-0.6)--(w);
\node  at  (-0.62,-0.2) {$x$};
\node   at   (0.15,0.1) {$\theta$};
\node  at  (0.24,-0.18) {$r$};
\draw [rotate=45]  [->] (0,0)--(w);
\node at (-0.2,-0.38) {$\sqrt{1+r^2}$};
\draw [->] (-1,0)--(1,0);
\draw [->] (0,-1)--(0,1);
\node at (-0.1,1) {$x_3$};
\node at (1,-0.1) {$x_1$};
\end{tikzpicture}
\caption{A spindle torus with axis of rotation $\theta$, tube centre offset $r$ and tube radius $\sqrt{1+r^2}$. The distance between the origin and either of the points where the torus self intersects is 1.}
\label{fig1}
\end{figure}
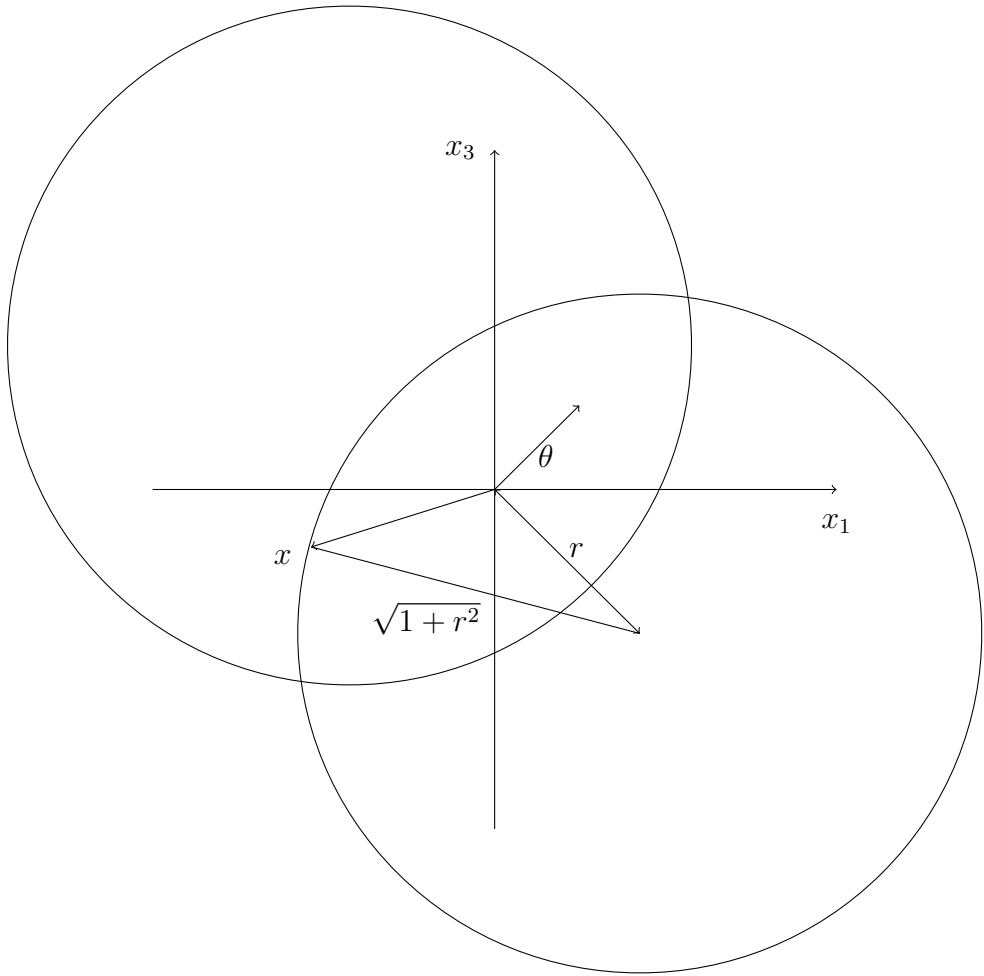

Here we present a microlocal analysis of the spindle transform, first introduced by the authors in \cite{WL}, which describes the Compton scattering tomography problem in three dimensions for a monochromatic source and energy sensitive detector pair. Compton scattering is the process in which a photon interacts in an inelastic collision with a charged particle. As the collision is inelastic, the photon undergoes a loss in energy, described by the equation
\begin{equation}
\label{equ1}
E_s=\frac{E_{\lambda}}{1+\left(E_{\lambda}/E_0\right)\left(1-\cos\omega\right)},
\end{equation}
where $E_s$ is the energy of the scattered photon which had an initial energy $E_{\lambda}$, $\omega$ is the scattering angle and $E_0\approx 511$keV is the electron rest energy. For $E_s$ and $E_{\lambda}$ fixed (i.e if the source is monochromatic and we can measure $E_s$), the scattering angle $\omega$ remains fixed and, in three dimensions, the surface of scatterers is the surface of revolution of a circular arc \cite{WL}. The surface of revolution of a circular arc is a \emph{spindle torus}. We define
\begin{equation}
T_{r}=\left\{(x_1,x_2,x_3)\in \mathbb{R}^3 : \left(r-\sqrt{x_1^2+x_2^2}\right)^2+x_3^2=1+r^2\right\}
\end{equation}
to be the spindle torus, radially symmetric about the $x_3$ axis, with tube centre offset $r\geq 0$ and tube radius $\sqrt{1+r^2}$. See figure \ref{fig1} which displays a rotation of $T_r$. 

In \cite{norton} Norton considered the problem of reconstructing a density supported in a quadrant of the plane from the Compton scattered intensity measured at a single point detector moved laterally along the axis away from a point source at the origin. Here the curve of scatterers is a circle. He considers the circle transform
\begin{equation}
Af(r,\phi)=\int_{-\frac{\pi}{2}}^{\frac{\pi}{2}}rF(r\cos\varphi,\varphi+\phi)\mathrm{d}\varphi,
\end{equation}
where $F(\rho,\theta)=f(\rho\cos\theta,\rho\sin\theta)$ is the polar form of $f :\mathbb{R}^2\to \mathbb{R}$. Let $f(x)=\frac{1}{|x|^2}\tilde{f}(\frac{x}{|x|^2})$. Then $Af(r,\phi)=R\tilde{f}(\frac{1}{r},\phi)$, where $R$ denotes the polar form of the straight line Radon transform. So $A$ is equivalent to $R$ via the diffeomorphism $x\to \frac{x}{|x|^2}$ and from this we can derive stability estimates from known theory on the Radon transform \cite{nat1}.  

In \cite{pal1,NT}, Nguyen and Truong consider an acquisition geometry of a point source and detector which remain opposite one another and are rotated on $S^1$, and aim to reconstruct a density supported on the unit disc. Here the curve of scatterers is a circular arc. They define the circular arc transform
\begin{equation}
Bf(r,\phi)=\int_{-\frac{\pi}{2}}^{\frac{\pi}{2}}\rho\sqrt{\frac{1+r^2}{1+r^2\cos^2{\varphi}}}F(\rho,\varphi+\phi)\mid_{\rho=\sqrt{r^2\cos^2 \varphi +1}-r\cos \varphi}\mathrm{d}\varphi.
\end{equation}
Letting $\tilde{f}(x)=\frac{\sqrt{|x|^2+1}-|x|}{1-|x|^2}f\left(\left(\sqrt{|x|^2+1}-|x|\right)\cdot \frac{x}{|x|}\right)$, we have that $\frac{r}{\sqrt{1+r^2}}Bf(r,\phi)=A\tilde{f}(r,\phi)$, and hence $B$ is equivalent to $A$ via the diffeomorphism $x\to \left(\sqrt{|x|^2+1}-|x|\right)\frac{x}{|x|}$. So, in two dimensions, the inverse problem for Compton scattering tomography is injective and mildly ill posed (the solution is bounded in some Sobolev space).  In \cite{pal1}, Palamodov also derives stability estimates for a more general class of Minkowski--Funk transforms.

In \cite{WL} the authors consider a three dimensional acquisition geometry, where a single source and detector are rotated opposite one another on $S^2$ and a density supported on a hollow ball is to be recovered. Here the surface of scatterers is a spindle torus. They define the spindle transform
\begin{equation}
\mathcal{S}f(r,\theta)=\int_{0}^{2\pi}\int_{0}^{\pi}\rho^2\sin\varphi\sqrt{\frac{1+r^2}{1+r^2\sin^2{\varphi}}} (h\cdot F)\left(\rho,\psi,\varphi\right)\mid_{\rho=\sqrt{r^2\sin^2 \varphi +1}-r\sin \varphi}\mathrm{d}\varphi\mathrm{d}\psi,
\end{equation}
where $F(\rho,\psi,\varphi)=f(\rho\cos\psi\sin\varphi,\rho\sin\psi\sin\varphi,\rho\cos\varphi)$ is the spherical polar form of $f :\mathbb{R}^3\to \mathbb{R}$ and $h\in\text{SO}(3)$ describes the rotation of the north pole to $\theta$, where $h$ defines a group action on real-valued functions in the natural way $(h\cdot f)(x)=f(hx)$. They show that a left inverse to $\mathcal{S}$ exists through the explicit inversion of a set of one-dimensional Volterra integral operators, and show that the null space of $\mathcal{S}$ consists of those functions whose even harmonic components are zero (odd functions). However the stability of the spindle transform was not considered. We aim to address this here from a microlocal perspective. In \cite{felea}, various acquisition geometries are considered for synthetic aperture radar imaging of moving objects. In each case the microlocal properties of the forward operator in question and its normal operator are analysed. In three of the four cases considered the Schwartz kernel of the normal operator was shown to belong to a class of distributions associated to two cleanly intersecting Lagrangians $I^{p,l}(\Delta, \Lambda)$ (that is, the wavefront set of the kernel of the normal operator is contained in $\Delta\cup\Lambda$). We show a similar result for the spindle transform $\mathcal{S}^*\mathcal{S}$, although the diagonal $\Delta$ is replaced by the disjoint union $\Delta \cup \widetilde{\Delta}$ where $\widetilde{\Delta}$ is reflection through the origin. We also determine the associated Lagrangian $\Lambda$. In \cite{felea} they suggest a way to reduce the size of the image artefact microlocally by applying an appropriate pseudodifferential operator as a filter before applying the backprojection operator. Similarly we derive a suitable filter for the spindle transform and show how it can be applied using the spherical harmonics of the data.

In section \ref{FIOsection} we show that $\mathcal{S}$ is equivalent to a weighted cylinder transform $\mathcal{C}$, which gives the weighted integrals over cylinders with an axis of revolution through the origin. After this we prove that $\mathcal{C}$ is a Fourier integral operator and determine its canonical relation. Later in section \ref{pairlagrangian} we present our main theorem (Theorem \ref{maintheorem}), where we show that $\mathcal{C}$ belongs to a class of distributions $I^{p,l}(\Delta\cup\widetilde{\Delta},\Lambda)$, where $\widetilde{\Delta}$ is a reflection and the Lagrangian $\Lambda$ is associated to a rotation artefact.

In section \ref{filtersection}, we adopt the ideas of Felea et al in \cite{felea} and derive a suitable pseudodifferential operator $Q$ which, when applied as a filter before applying the backprojection operator of the cylinder transform, reduces the artefact intensity in the image. We show that $Q$ can be applied by multiplying the harmonic components of the data by a factor $c_l$, which depends on the degree $l$ of the component, and show how this translates to a spherical convolution of the data with a distribution on the sphere $h$.

Simulated reconstructions from spindle transform data are presented in section \ref{results}. We reconstruct a small bead of constant density by unfiltered backprojection and show the artefacts produced by $\Lambda$ in our reconstruction. We then reconstruct the same density by filtered backprojection, applying the filter $Q$ as an intermediate step, and show how the size of the artefacts are reduced in the image. Later we provide reconstructions of densities of oscillating layers using the conjugate gradient least squares (CGLS) method and Landweber iteration. We arrange the layers as spherical shells centred at the origin and as planes and  compare our results. We also investigate the effects of applying the filter $Q$ as a pre--conditioner, prior to implementing CGLS and the Landweber method. 

\section{The microlocal properties of $\mathcal{S}$ and $\mathcal{S}^*\mathcal{S}$}
\label{mainsection}
Here we investigate the microlocal properties of the spindle transform and its normal operator. We start by showing the equivalence of $\mathcal{S}$ to a cylinder transform $\mathcal{C}$, and how we can write $\mathcal{S}$ and $\mathcal{C}$ as Fourier integral operators. Then we determine the canonical relations associated with $\mathcal{C}$ and from these we discover that $\mathcal{C}^*\mathcal{C}$ is a paired Lagrangian operator with blowdown--blowdown singularities. First we give some preliminaries.

Let $B^n_{\epsilon_1,\epsilon_2}=\{x\in \mathbb{R}^n : 0<\epsilon_1 <|x|<\epsilon_2<1\}$ denote the set of points on a hollow ball with inner radius $\epsilon_1$ and outer radius $\epsilon_2$. Let $Z^n=\mathbb{R}\times S^{n-1}$ denote the $n$--cylinder and, for $X\subset \mathbb{R}^n$ an open set, let $\mathcal{D}'(X)$ denote the vector space of distributions on $X$, and let $\mathcal{E}'(X)$ denote the vector space of distributions with compact support contained in $X$.
\begin{definition}
For a function $f$ in the Schwarz space $S(\mathbb{R}^n)$ we define the Fourier transform and its inverse in terms of angular frequency as
\begin{equation}
\begin{split}
\mathcal{F}f(\xi)&=(2\pi)^{-\frac{n}{2}}\int_{\mathbb{R}^n}e^{-ix\cdot\xi}f(x)\mathrm{d}x,\\
\mathcal{F}^{-1}f(x)&=(2\pi)^{-\frac{n}{2}}\int_{\mathbb{R}^n}e^{ix\cdot\xi}f(\xi)\mathrm{d}\xi.
\end{split}
\end{equation}
\end{definition}
\begin{definition}
Let $m,\rho,\delta \in\mathbb{R}$ with $0\leq\rho\leq1$ and $\delta=1-\rho$. Then we define $S^m_{\rho}(X\times\mathbb{R}^n)$ to be the set of $a\in C^{\infty}(X\times \mathbb{R}^n)$ such that for every compact set $K\subset X$ and all multi--indices $\alpha, \beta$ the bound
\begin{equation}
\left|\partial^{\beta}_x\partial^{\alpha}_{\xi}a(x,\xi)\right|\leq C_{\alpha,\beta,K}(1+|\xi|)^{m-\rho|\alpha|+\delta|\beta|},\ \ \ x\in K,\ \xi\in\mathbb{R}^n,
\end{equation}
holds for some constant $C_{\alpha,\beta,K}$. The elements of $S^m_{\rho}$ are called \emph{symbols} of order $m$, type $\rho$.
\end{definition}
\begin{definition}
A function $\phi=\phi(x,\xi)\in C^{\infty}(X\times\mathbb{R}^N\backslash 0)$ is a \emph{phase function} if $\phi(x,\lambda\xi)=\lambda\phi(x,\xi)$, $\forall \lambda>0$ and $\mathrm{d}\phi\neq 0$.
\end{definition}
\begin{definition}
Let $X \subset \mathbb{R}^{n_x}$, $Y\in \subset \mathbb{R}^{n_y}$ be open sets. A Fourier integral operator (FIO) of order $m + N/2 - (n_x + n_y)/4$ is an operator $A:C^{\infty}_0(X)\to \mathcal{D}'(Y)$ with Schwartz kernel given by an oscillatory integral of the form
\begin{equation} \label{oscint}
Af(y)=\int_{\mathbb{R}^N} e^{i\phi(x,y,\xi)}a(x,y,\xi) \mathrm{d}\xi,
\end{equation}
where $\phi$ is a phase function, and $a \in S^m_{\rho}((X \times Y) \times \mathbb{R}^N)$ is a symbol.
\end{definition}
\begin{definition}\label{def:canon}
The canonical relation of an FIO with phase function $\phi$ is defined as
\begin{equation}
\begin{split}
C=&\{\left((y,\eta),(x,\omega)\right)\in\left(Y\times \mathbb{R}^N\backslash 0\right)\times\left(X\times\mathbb{R}^N\backslash 0\right) : (x,y,\omega)\in \Sigma_\phi,\\
&\omega=-d_x\phi(x,y,\xi), \eta=d_y\phi(x,y,\xi), \omega,\eta\neq0\},
\end{split}
\end{equation}
where $\Sigma_\phi=\{(x,\xi)\in X\times\mathbb{R}^N\backslash 0 : d_{\xi}\phi=0\}$ is the critical set of $\phi$.
\end{definition}

\noindent If $Y$ and $X$ are manifolds without boundary, then an operator $A: C_0^\infty(X) \rightarrow \mathcal{D}'(Y)$ is a Fourier integral operator if its Schwartz kernel can be represented locally in coordinates by oscillatory integrals of the form \eqref{oscint}, and the canonical relations of the phase functions for the local representations all lie within a single immersed Lagrangian submanifold of $T^*Y \times T^*X$. For much more detail on Fourier integral operators and their definition see \cite{hormander}.
\subsection{The spindle transform as an FIO}
\label{FIOsection}
Recall the author's acquisition geometry in \cite{WL} (displayed in figure \ref{fig1}). We have the implicit equation
\begin{equation}
(r+|x\times\theta|)^2+(x\cdot\theta)^2=1+r^2
\end{equation}
for the set of points on a spindle with tube centre offset $r$ and axis of revolution given by $\theta\in S^2$. With this in mind we define
\begin{equation}
h(s,x,\theta)=\frac{4|x\times\theta|^2}{(1-|x|^2)^2}-s,
\end{equation}
and then we can write the spindle transform $\mathcal{S}:C^{\infty}_0(B^3_{\epsilon_1,\epsilon_2})\to C^{\infty}((0,1)\times S^2)$ as
\begin{equation}
\begin{split}
\label{spindle}
\mathcal{S}f(s,\theta)&=\int_{B^3_{\epsilon_1,\epsilon_2}}\frac{\delta\left(\frac{4|x\times\theta|^2}{(1-|x|^2)^2}-s\right)}{\left|\nabla_x h(s,x,\theta)\right|}f(x)\mathrm{d}x,
\end{split}
\end{equation}
where $s=1/r^2$ and $\delta$ is the Dirac--delta function. Note that
\begin{equation} \label{nabh}
\nabla_x h(s,x,\theta) = 8\frac{(x-(x\cdot\theta)\theta)}{(1-|x|^2)^2}+16 \frac{|x\times\theta|^2x}{(1-|x|^2)^3}
\end{equation}
is smooth, bounded above, and does not vanish on on $(0,1) \times B^3_{\epsilon_1,\epsilon_2} \times S^2$. We define the backprojection operator $\mathcal{S}^*:C^{\infty}((0,1)\times S^2)\to C^{\infty}(B^3_{\epsilon_1,\epsilon_2})$ as
\begin{equation}
\begin{split}
\mathcal{S}^*g(x)&=\int_{S^2}\frac{g\left(\frac{4|x\times\theta|^2}{(1-|x|^2)^2},\theta\right)}{\left|\nabla_x h(s,x,\theta)\right|}\mathrm{d}\Omega,
\end{split}
\end{equation}
where $\mathrm{d}\Omega$ is the surface measure on $S^2$.

\begin{proposition}
\label{prop1}
The backprojection operator $\mathcal{S}^*$ is the adjoint operator to $\mathcal{S}$.
\begin{proof}
Let $g\in C^{\infty}((0,1)\times S^2)$ and $f\in C^{\infty}_0(B^3_{\epsilon_1,\epsilon_2})$. Then (in the third step note that $\nabla_x h$ does not actually depend on $s$)
\begin{equation}
\begin{split}
\langle g, \mathcal{S}f \rangle &=\int_{S^2}\int_0^1 g(s,\theta)\mathcal{S}f(s,\theta)\mathrm{d}s\mathrm{d}\Omega\\
&=\int_{S^2}\int_0^1 g(s,\theta)\int_{B^3_{\epsilon_1,\epsilon_2}}\frac{\delta\left(\frac{4|x\times\theta|^2}{(1-|x|^2)^2}-s\right)}{\left|\nabla_x h(s,x,\theta)\right|}f(x)\mathrm{d}x\mathrm{d}s\mathrm{d}\Omega\\
&=\int_{B^3_{\epsilon_1,\epsilon_2}}\int_{S^2}\frac{g\left(\frac{4|x\times\theta|^2}{(1-|x|^2)^2},\theta\right)}{\left|\nabla_x h(s,x,\theta)\right|}\mathrm{d}\Omega f(x)\mathrm{d}x\\
&=\int_{B^3_{\epsilon_1,\epsilon_2}}\mathcal{S}^*g(x)f(x)\mathrm{d}x=\langle \mathcal{S}^*g, f\rangle,
\end{split}
\end{equation}
which completes the proof.
\end{proof}
\end{proposition}

Let $v(x)=\left(\sqrt{1+\frac{1}{|x|^2}}-\frac{1}{|x|}\right)\cdot \frac{x}{|x|}$, and set $\alpha_i = 2 \epsilon_i/(1-\epsilon_i^2)$ for $i = 1$ or $2$ so that when $|x| = \alpha_i$, $|v(x)| = \epsilon_i$. Then, after making the substitution $x\to v(x)$ in equation (\ref{spindle}), we have
\begin{equation}
\label{cylinder}
\begin{split}
\mathcal{S}f(s,\theta)&=\int_{B^3_{\alpha_1,\alpha_2}}|\text{det}(J_v)|\frac{\delta\left(|x\times \theta|^2-s\right)}{|\nabla_v h(s,v(x),\theta)|}f\left(\left(\sqrt{1+\frac{1}{|x|^2}}-\frac{1}{|x|}\right)\cdot \frac{x}{|x|}\right)\mathrm{d}x\\
&=\int_{B^3_{\alpha_1,\alpha_2}}\frac{\delta\left(|x|^2-(x\cdot\theta)^2-s\right)}{|\nabla_v h(s,v(x),\theta)|}\tilde{f}(x)\mathrm{d}x,\\
\end{split}
\end{equation}
where
\begin{equation}
\tilde{f}(x)=|\text{det}(J_v)|\ f\left(\left(\sqrt{1+\frac{1}{|x|^2}}-\frac{1}{|x|}\right)\cdot \frac{x}{|x|}\right). 
\end{equation}
We define the weighted cylinder transform $\mathcal{C}:C^{\infty}_0(B^3_{\alpha_1,\alpha_2})\to C^{\infty}((0,1)\times S^2)$ as
\begin{equation}
\label{cylinder1}
\mathcal{C}f(s,\theta)=\int_{B^3_{\alpha_1,\alpha_2}}\frac{\delta\left(|x|^2-(x\cdot\theta)^2-s\right)}{\sqrt{|\nabla_v h(s,v(x),\theta)|}}f(x)\mathrm{d}x
\end{equation}
and its backprojection operator $\mathcal{C}^*:C^{\infty}((0,1)\times S^2)\to C^{\infty}(B^3_{\alpha_1,\alpha_2})$:
\begin{equation}
\mathcal{C}^*g(x)=\int_{S^2}\frac{g\left(|x\times\theta|^2,\theta\right)}{|\nabla_v h(s,v(x),\theta)|}\mathrm{d}\Omega.
\end{equation}
As in Proposition \ref{prop1}, we can show that $\mathcal{C}^*$ is the formal adjoint to $\mathcal{C}$.

The above is to say that the spindle transform is equivalent, via the diffeomorphism $x\to\left(\sqrt{1+\frac{1}{|x|^2}}-\frac{1}{|x|}\right)\cdot \frac{x}{|x|}$, to the transform $\mathcal{C}$ which defines the weighted integrals over cylinders with radius $\sqrt{s}$ and axis of rotation through the origin with direction $\theta$. With this in mind we consider the microlocal properties of the cylinder transform $\mathcal{C}$ and its normal operator $\mathcal{C}^*\mathcal{C}$ for the remainder of this section.

First, we characterise $\mathcal{C}$ as a Fourier integral operator in the next theorem.
\begin{theorem}
\label{FIO}
The cylinder transform $\mathcal{C}$ is a Fourier integral operator order $-1$ with canonical relation
\begin{equation}
\begin{split}
C=&\big\{\left((s,\alpha,\beta),(\sigma,2\sigma(x\cdot\theta_{\alpha})(x\cdot\theta),2\sigma(x\cdot\theta_{\beta})(x\cdot\theta));x,2\sigma(x-(x\cdot\theta)\theta)\right) : x\in B^3_{\alpha_1,\alpha_2},\\
&s\in(0,1),\sigma\in\mathbb{R}\backslash0,\theta\in S^2, |x|^2-(x\cdot\theta)^2-s=0\big\},
\end{split}
\end{equation}
where $(\alpha, \beta) \in \mathbb{R}^2$ provide a local parameterization of $\theta$, $\theta_{\alpha}=\partial_{\alpha}\theta$ and $\theta_{\beta}=\partial_{\beta}\theta$.
\begin{proof}
The delta function may be written as the oscillatory integral
\begin{equation}
\delta(s)=\frac{1}{2\pi}\int_{-\infty}^{\infty}e^{i\sigma s}\mathrm{d}\sigma.
\end{equation}
Thus, by equation (\ref{cylinder1}) we have
\begin{equation}
\begin{split}
\mathcal{C}f(s,\theta)&=\int_{B^3_{\alpha_1,\alpha_2}}\frac{\delta\left(|x|^2-(x\cdot\theta)^2-s\right)}{|\nabla_v h(s,v(x),\theta)|}f(x)\ \mathrm{d}x\\
&=\frac{1}{2\pi}\int_{-\infty}^{\infty}\int_{B^3_{\alpha_1,\alpha_2}}\frac{e^{i(s-|x|^2+(x\cdot\theta)^2)\sigma}}{|\nabla_v h(s,v(x),\theta)|}f(x)\ \mathrm{d}x\mathrm{d}\sigma,
\end{split}
\end{equation}
and from this we see that $\mathcal{C}$ is an FIO with phase function
\begin{equation}
\phi(x,s,\theta,\sigma)=(s-|x|^2+(x\cdot\theta)^2)\sigma,
\end{equation}
and amplitude
\begin{equation}
a(x,s,\theta,\sigma)= \frac{1}{2\pi} |\nabla_v h(s,v(x),\theta)|^{-1}
\end{equation}
where the single phase variable is $\sigma$. Indeed, as we noted above and is evident from the formula \eqref{nabh}, $|\nabla_v h(s,v,\theta)|$ is smooth, bounded from above, and bounded from below larger than zero when $v \in B^3_{\epsilon_1,\epsilon_2}$. Also, $x \mapsto v(x)$ is a diffeomorphism from $B^3_{\alpha_1,\alpha_2}$ to $B^3_{\epsilon_1,\epsilon_2}$. Therefore, since $a$ also does not depend on the phase variable $\sigma$, $a\in S^0_{1}((B_{\alpha_1,\alpha_2}\times (0,1)\times S^2) \times \mathbb{R})$. Hence the order of $\mathcal{C}$ is $0+\frac{1}{2}-\frac{1}{4}(3+3)=-1$.

Now suppose that $\theta \in S^2$ is parametrized by $\alpha$ and $\beta \in \mathbb{R}$ (for example using standard spherical coordinates). Then
\begin{equation}
\phi(x,s,\alpha,\beta,\sigma)=(s-|x|^2+(x\cdot\theta)^2)\sigma,
\end{equation}
and the derivatives of $\phi$ are
\begin{equation}
\begin{split}
\mathrm{d}_x\phi& =-2\sigma(x-(x\cdot\theta)\theta),\ \ \mathrm{d}_{\alpha}\phi=2\sigma(x\cdot\theta_{\alpha})(x\cdot\theta),\ \ \mathrm{d}_{\beta}\phi=2\sigma(x\cdot\theta_{\beta})(x\cdot\theta),\\&
\mathrm{d}_s\phi=\sigma, \ \ \mathrm{d}_{\sigma}\phi=s-|x|^2+(x\cdot\theta)^2.
\end{split}
\end{equation}
From Definition \ref{def:canon}, it follows that the canonical relation of $\mathcal{C}$ is:
\begin{equation}
\hspace*{-1cm}
\begin{split}
C&=\Big\{\left((s,\alpha,\beta),(\hat{s},\hat{\alpha},\hat{\beta}) ; (x,\xi)\right) : x\in B^3_{\alpha_1,\alpha_2},\\
&\ \ \  \ \ \  \ \ \ \ \ \  s\in(0,1),\sigma\in\mathbb{R}\backslash0,\theta\in S^2, \xi=-\mathrm{d}_x\phi, \hat{s}=\mathrm{d}_s\phi, \hat{\alpha}=\mathrm{d}_{\alpha}\phi, \hat{\beta}=\mathrm{d}_{\beta}\phi, \mathrm{d}_{\sigma}\phi=0\Big\}\\
&=\Big\{\left( (s,\alpha,\beta),(\sigma,2\sigma(x\cdot\theta_{\alpha})(x\cdot\theta),2\sigma(x\cdot\theta_{\beta})(x\cdot\theta)) ;x,2\sigma(x-(x\cdot\theta)\theta)\right) : x\in B^3_{\alpha_1,\alpha_2},\\
&\ \ \  \ \ \  \ \ \ \ \ \ s\in(0,1),\sigma\in\mathbb{R}\backslash0,\theta\in S^2, |x|^2-(x\cdot\theta)^2-s=0\Big\},
\end{split}
\end{equation}
which completes the proof.
\end{proof}
\end{theorem}
\subsection{$\mathcal{C}^*\mathcal{C}$ as a paired Lagrangian operator}
\label{pairlagrangian}
If we analyse the canonical relation $C$ given in Theorem \ref{FIO}, we can see that it is non-injective as the points on the ring $\{x\in \mathbb{R}^3 :  |x|^2-s=0, x\cdot \theta=0\}$ map to $\left((s,\theta),(\sigma,0)\right)$ if we fix $s$ and $\sigma$. Let $C^*$ be the canonical relation of $\mathcal{C}^*$ and let $\Delta$ denote the diagonal. Then, given the non-injectivity of $C$, $C^*\circ C\nsubseteq \Delta$ and $\mathcal{C}^*\mathcal{C}$ is not a pseudodifferential operator, or even an FIO. In this section we show that the Schwarz kernel of $\mathcal{C}^*\mathcal{C}$ instead belongs to a class of distributions $I^{p,l}(\Delta\cup\widetilde{\Delta}, \Lambda)$ studied in \cite{felea,theoremref}. First we recall some definitions and theorems from \cite{felea}.
\begin{definition}
Two submanifolds $M, N\subset X$ intersect cleanly if $M\cap N$ is a smooth submanifold and $T(M\cap N)=TM\cap TN$
\end{definition}
\begin{definition}
We define $I^m(C)$ to be the set of Fourier integral operators, $A : \mathcal{E}'(X)\to \mathcal{D}'(Y)$, of order $m$ with canonical relation $C\subset( T^*Y\setminus 0) \times (T^*X \setminus 0)$
\end{definition}
\noindent Recall the definitions of the left and right projections of a canonical relation.
\begin{definition}
Let $C$ be the canonical relation associated to the FIO $A:\mathcal{E}'(X)\to \mathcal{D}'(Y)$. Then we denote $\pi_L$ and $\pi_R$ to be the left and right projections of $C$, $\pi_L:C\to T^*Y\backslash 0$ and $\pi_R : C\to T^*X\backslash 0$.
\end{definition}
\noindent We have the following result from \cite{hormander}.
\begin{proposition}
Let $\text{dim}(X)=\text{dim}(Y)$. Then at any point in $C$:
\begin{enumerate}[(i)]
\item if one of $\pi_L$ or $\pi_R$ is a local diffeomorphism, then $C$ is a local canonical graph;
\item if one of the projections $\pi_R$ or $\pi_L$ is singular, then so is the other. The type of the singularity may be different (e.g. fold or blowdown \cite{felea1}) but both projections drop rank on the same set
\begin{equation}
\Sigma=\{(y,\eta; x,\xi)\in C : \text{det}(\mathrm{d}\pi_L)=0\}=\{(y,\eta; x,\xi)\in C : \text{det}(\mathrm{d}\pi_R)=0\}.
\end{equation}
\end{enumerate}
\end{proposition}
\noindent Now we have the definition of a blowdown singularity and the definitions of a nonradial and involutive submanifold:
\begin{definition}
Let $M$ and $N$ be manifolds of dimension $n$ and let $f:N\to M$ be a smooth function. $f$ is said to have a \emph{blowdown} singularity of order $k \in \mathbb{N}$ along a smooth hypersurface $\Sigma\subset M$ if $f$ is a local diffeomorphism away from $\Sigma$, $\mathrm{d}f$ drops rank by $k$ at $\Sigma$, $\text{ker}(\mathrm{d}f)\subset T(\Sigma)$, and the determinant of the Jacobian matrix vanishes to order $k$ at $\Sigma$.
\end{definition}
\begin{definition}
A submanifold $M\subset T^*X$ is \emph{nonradial} if $\rho\notin (TM)^{\perp}$, where $\rho=\sum\xi_i\partial_{\xi_i}$.
\end{definition}
\begin{definition}
A submanifold $M\subset T^*X$, $M=\{(x,\xi) : p_i(x,\xi)=0, 1\leq i\leq k\}$ is \emph{involutive} if the differentials $\mathrm{d}p_i$, $i=1,\ldots, k$, are linearly independant and the Poisson brackets satisfy $\{p_i,p_j\}=0$, $i\neq j$.
\end{definition}
\noindent From \cite{greenleaf}, we have the definition of the flowout.
\begin{definition}
Let $\Gamma=\{(x,\xi) : p_i(x,\xi)=0, 1\leq i\leq k\}$ be a submanifold of $T^*X$. Then the \emph{flowout} of $\Gamma$ is given by $\{(x,\xi;y,\eta)\in T^*X\times T^*X : (x,\xi)\in \Gamma, (y,\eta)=\exp(\sum_{i=1}^kt_iH_{p_i})(x,\xi), t\in\mathbb{R}^k\}$, where $H_{p_i}$ is the Hamiltonian vector field of $p_i$.
\end{definition}
\noindent We now state a result of \cite[Theorem 1.2]{theoremref} concerning the composition of FIO's with blowdown--blowdown singularities.
\begin{theorem}
\label{blowdowntheorem}
Let $C\subset (T^*Y\backslash 0)\times (T^*X\backslash 0)$ be a canonical relation which satifies the following:
\begin{enumerate}[(i)]
\item away from a hypersurface $\Sigma\subset C$, the left and right projections $\pi_L$ and $\pi_R$ are diffeomorphisms;
\item at $\Sigma$, both $\pi_L$ and $\pi_R$ have blowdown singularities;
\item $\pi_L(\Sigma)$ and $\pi_R(\Sigma)$ are nonradial and involutive.
\end{enumerate}
If $A\in I^m(C)$ and $B\in I^{m'}(C^t)$, then $BA\in I^{m+m'+\frac{k-1}{2},-\frac{k-1}{2}}(\Delta,\Lambda_{\pi_R(\Sigma)})$, where $\Delta$ is the diagonal and $\Lambda_{\pi_R(\Sigma)}$ is the flowout of of $\pi_R(\Sigma)$.
\end{theorem}
\noindent Finally, for  two cleanly intersecting Lagrangians $\Lambda_0$ and $\Lambda_1$, we define the $I^{p,l}(\Lambda_0,\Lambda_1)$ classes as in  \cite{felea,IP}.

We now have our main Theorem.
\begin{theorem}
\label{maintheorem}
Let $C$ be the canonical relation of the cylinder transform $\mathcal{C}$. Then the left and right projections of $C$ have singularities along a dimension 1 submanifold $\Sigma$, $\pi_L(\Sigma)$ and $\pi_R(\Sigma)$ are involutive and nonradial, and $\mathcal{C}^*\mathcal{C}\in I^{-2,0}(\Delta\cup\widetilde{\Delta},\Lambda)$, where $\Delta$ is the diagonal in $T^*B^3_{\alpha_1,\alpha_2} \times T^*B^3_{\alpha_1,\alpha_2}$, $\widetilde{\Delta} = \{(-x,-\xi; x, \xi) \ : \ (x,\xi) \in T^* B^3_{\alpha_1,\alpha_2}\}$, and $\Lambda$ is the flowout of $\pi_R(\Sigma)$.
\begin{proof}
From Theorem \ref{FIO} we have the canonical relation of the cylinder transform
\begin{equation}
\begin{split}
C=&\big\{\left( (s,\alpha,\beta),(\sigma,2\sigma(x\cdot\theta_{\alpha})(x\cdot\theta),2\sigma(x\cdot\theta_{\beta})(x\cdot\theta)); x,2\sigma(x-(x\cdot\theta)\theta) \right) : x\in B^3_{\alpha_1,\alpha_2},\\
&s\in (0,1),\sigma\in\mathbb{R}\backslash0,\theta\in S^2, |x|^2-(x\cdot\theta)^2-s=0\big\},
\end{split}
\end{equation}
where $\alpha$ and $\beta$ parameterize $\theta \in S^2$.
Suppose we use standard spherical coordinates centred at any given point on $S^2$ (e.g. when centred at $(1,0,0)$ these would be defined by $\theta = (\cos\alpha\cos\beta,\sin\alpha\cos\beta,\sin\beta)$) with the notation
\begin{equation}
\partial_{\alpha}\theta=\theta_{\alpha},\ \ \ \partial_{\beta}\theta=\theta_{\beta}.
\end{equation}
With such a parameterization, $\{\theta,\theta_{\alpha},\theta_{\beta}\}$ is an orthogonal basis for $\mathbb{R}^3$ and $|\theta_{\beta}|=1$, $|\theta_{\alpha}|=\cos\beta$. Furthermore
\begin{equation}
\text{det}\begin{pmatrix}
  \theta^{T}  \\
  \theta^{T}_{\alpha} \\
  \theta^{T}_{\beta}  
 \end{pmatrix}
=\cos\beta.
\end{equation}
Using these coordinates on $S^2$, we can also parametrise $C$ by $(x,\alpha,\beta,\sigma)$ where $x \in B_{\alpha_1,\alpha_2}^3$ is such that $|x|^2 - (x \cdot \theta)^2 = (\cos\beta)^{-2} (x\cdot \theta_\alpha)^2 + (x \cdot \theta_\beta)^2 \in (0,1)$, and $(\alpha,\beta)$ are in the domain of the coordinates for $S^2$. With this parametrization of $C$ the left projection is given by
\begin{equation}
\pi_L(x,\alpha,\beta,\sigma)=\left(|x|^2-(x\cdot \theta)^2,\alpha,\beta,\sigma,2\sigma (x\cdot\theta_{\alpha})(x\cdot\theta),2\sigma (x\cdot\theta_{\beta})(x\cdot\theta)\right).
\end{equation}
First note that $\pi_L(x,\alpha,\beta,\sigma) = \pi_L(-x,\alpha,\beta,\sigma)$, and so $\pi_L$ is not injective. However, as we shall see, except for on the set $\Sigma = \{x \cdot \theta\}$, $\pi_L$ is exactly two-to-one. Indeed, suppose that $x \cdot \theta \neq 0$ and $\pi_L(x,\alpha,\beta,\sigma) = \pi_L(-x',\alpha',\beta',\sigma')$. Then $\alpha = \alpha'$, $\beta = \beta'$, $\sigma = \sigma'$,
\[
|x|^2 - (x\cdot \theta)^2 = |x'|^2 - (x'\cdot \theta)^2 \Leftrightarrow (\cos \beta)^{-2} (x \cdot \theta_\alpha)^2 + (x \cdot \theta_\beta)^2 =(\cos \beta')^{-2} (x' \cdot \theta_\alpha)^2 + (x' \cdot \theta_\beta)^2,
\]
and (using the fact that $\sigma = \sigma' \neq 0$)
\[
(x \cdot \theta) (x\cdot \theta_\alpha, x \cdot \theta_\beta) = (x' \cdot \theta) (x'\cdot \theta_\alpha, x' \cdot \theta_\beta).
\]
Since $x \cdot \theta \neq 0$, and $(\cos \beta')^{-2} (x' \cdot \theta_\alpha)^2 + (x' \cdot \theta_\beta)^2 \neq 0$, we can combine these to see that $x = \pm x'$.

Now let us analyze $D \pi_L$ to show that $\pi_L$ is a local diffeomorphism away from $\Sigma$. Letting $I_{n\times n}$ and $0_{n\times n}$ denote the $n\times n$ identity and zero matrices respectively, after a permutation of rows, the differential of $\pi_L$ is
\begin{equation}
D \pi_L= \begin{pmatrix}
 2(x^T-(x\cdot\theta)\theta^T)  & r_1 \\
2\sigma((x\cdot\theta_{\alpha})\theta^T+(x\cdot\theta)\theta^T_{\alpha})  & r_2 \\
2\sigma((x\cdot\theta_{\beta})\theta^T+(x\cdot\theta)\theta^T_{\beta})  & r_3 \\
  0_{3\times 3}  & I_{3\times 3}
 \end{pmatrix},
\end{equation}
where
\begin{equation}
\begin{split}
r_1&=-\left(2(x\cdot \theta_{\alpha})(x\cdot\theta),2(x\cdot \theta_{\beta})(x\cdot\theta),0\right),\\
r_2&=\left(2\sigma((x\cdot \theta_{\alpha})^2+(x\cdot \theta_{\alpha\alpha})(x\cdot \theta)),2\sigma((x\cdot \theta_{\alpha})(x\cdot \theta_{\beta})+(x\cdot \theta_{\alpha\beta})(x\cdot \theta)),2(x\cdot\theta_{\alpha})(x\cdot\theta)\right),\\
r_3&=\left(2\sigma((x\cdot \theta_{\beta})(x\cdot \theta_{\alpha})+(x\cdot \theta_{\alpha\beta})(x\cdot \theta)),2\sigma((x\cdot \theta_{\beta})^2+(x\cdot \theta_{\beta\beta})(x\cdot \theta)),2(x\cdot\theta_{\beta})(x\cdot\theta)\right),
\end{split}
\end{equation}
and $\theta_{\alpha\alpha}=\partial_{\alpha\alpha}\theta$, $\theta_{\beta\beta}=\partial_{\beta\beta}\theta$ and $\theta_{\alpha\beta}=\partial_{\alpha\beta}\theta$.

We can now calculate the determinant of $D \pi_L$ as follows:
\begin{equation}
\begin{split}
\text{det}D \pi_L&=\text{det} \begin{pmatrix}
 2(x^T-(x\cdot\theta)\theta^T)  \\
2\sigma((x\cdot\theta_{\alpha})\theta^T+(x\cdot\theta)\theta^T_{\alpha})  \\
2\sigma((x\cdot\theta_{\beta})\theta^T+(x\cdot\theta)\theta^T_{\beta})
 \end{pmatrix}\\
&=\frac{1}{\cos\beta}\text{det}\begin{pmatrix}
 2(x^T-(x\cdot\theta)\theta^T)  \\
2\sigma((x\cdot\theta_{\alpha})\theta^T+(x\cdot\theta)\theta^T_{\alpha})  \\
2\sigma((x\cdot\theta_{\beta})\theta^T+(x\cdot\theta)\theta^T_{\beta})
 \end{pmatrix}\left(\theta,\theta_{\alpha},\theta_{\beta}\right)\\
&=\frac{1}{\cos\beta}\text{det}\begin{pmatrix}
 0 & 2(x\cdot\theta_{\alpha}) & 2(x\cdot\theta_{\beta})  \\
2\sigma(x\cdot\theta_{\alpha}) & 2\sigma(x\cdot\theta)\cos^2\beta & 0 \\
2\sigma(x\cdot\theta_{\beta}) & 0 & 2\sigma(x\cdot \theta)
\end{pmatrix}\\
&=\frac{8\sigma^2}{\cos\beta}(x\cdot\theta)\left((x\cdot\theta_{\alpha})^2+(x\cdot\theta_{\beta})^2\cos^2\beta\right).
\end{split}
\end{equation}
This is zero when $x\cdot\theta=0$ or $(x\cdot\theta_{\alpha})^2+(x\cdot\theta_{\beta})^2\cos^2\beta=0$. The latter case corresponds to when $x$ and $\theta$ are parallel, which we do not consider ($x$ and $\theta$ are parallel only when the cylinder is degenerate, i.e. when $s=0$), and so $\pi_L$ is a local diffeomorphism away from the manifold $\Sigma=\{x\cdot\theta=0\}$.

Finally we show that the singularities of $\pi_L$ on $\Sigma$ are blowdown of order $1$. Indeed, on $\Sigma$ we have
\begin{equation}
\mathrm{d}\ \text{det} D\pi_L=\frac{8\sigma^2}{\cos\beta}\left((x\cdot\theta_{\alpha})^2+(x\cdot\theta_{\beta})^2\cos^2\beta\right)\left(\theta\cdot\mathrm{d}x+(x\cdot\theta_{\alpha})\mathrm{d}\alpha+(x\cdot\theta_{\beta})\mathrm{d}\beta\right),
\end{equation}
and the kernel of $D\pi_L$ on $\Sigma$ is
\begin{equation}
\text{span}\left\{\left((x\cdot\theta_{\beta})\theta_{\alpha}-(x\cdot\theta_{\alpha})\theta_{\beta}\right)\cdot\nabla_x\right\}\subset \ker\left\{\mathrm{d}\ \text{det} D\pi_L\right\}.
\end{equation}
So the left projection $\pi_L$ drops rank by 1 on $\Sigma$ and its critical points on $\Sigma$ are blowdown type singularities. Furthermore,
\begin{equation}
\pi_L(\Sigma)=\{\hat{\alpha}=\hat{\beta}=0\},
\end{equation}
which is involutive and nonradial (here $\hat{\alpha}$ and $\hat{\beta}$ are the dual variables of $\alpha$ and $\beta$).

Using the same parameterization of $C$ as above, the right projection is given by
\begin{equation}
\pi_R(x,\alpha,\beta,\sigma)=\left(x,2\sigma(x-(x\cdot\theta)\theta)\right),
\end{equation}
and its differential is
\begin{equation}
\hspace*{-0.75cm}
D\pi_R=\begin{pmatrix}
I_{3\times 3} & 0_{3\times 1} & 0_{3\times 1} & 0_{3\times 1}\\
2\sigma\left(I_{3\times 3}-\theta\theta^T\right) & -2\sigma((x\cdot\theta_{\alpha})\theta+(x\cdot\theta)\theta_{\alpha}) & -2\sigma((x\cdot\theta_{\beta})\theta+(x\cdot\theta)\theta_{\beta}) & 2(x-(x\cdot\theta)\theta).
\end{pmatrix}.
\end{equation}
The determinant can now be calculated as
\begin{equation}
\hspace*{-0.75cm}
\begin{split}
\text{det} D\pi_R&=\text{det} \left(-2\sigma((x\cdot\theta_{\alpha})\theta+(x\cdot\theta)\theta_{\alpha}),-2\sigma((x\cdot\theta_{\beta})\theta+(x\cdot\theta)\theta_{\beta}),2(x-(x\cdot\theta)\theta)\right)\\
&=\frac{1}{\cos\beta}\text{det}\begin{pmatrix}
  \theta^{T}  \\
  \theta^{T}_{\alpha} \\
  \theta^{T}_{\beta}  
 \end{pmatrix}\left(-2\sigma((x\cdot\theta_{\alpha})\theta+(x\cdot\theta)\theta_{\alpha}),-2\sigma((x\cdot\theta_{\beta})\theta+(x\cdot\theta)\theta_{\beta}),2(x-(x\cdot\theta)\theta)\right)\\
&=\frac{1}{\cos\beta}\text{det}\begin{pmatrix}
-2\sigma(x\cdot\theta_{\alpha}) & -2\sigma(x\cdot \theta_{\beta}) & 0\\
-2\sigma(x\cdot\theta)\cos^2\beta & 0 & 2(x\cdot\theta_{\alpha})\\
0 & -2\sigma(x\cdot\theta) & 2(x\cdot\theta_{\beta})
\end{pmatrix}\\
&=-\frac{8\sigma^2}{\cos\beta}(x\cdot\theta)\left((x\cdot\theta_{\alpha})^2+(x\cdot\theta_{\beta})^2\cos^2\beta\right).
\end{split}
\end{equation}
Hence, on $\Sigma$
\begin{equation}
\mathrm{d}\ \text{det} D\pi_R=\frac{8\sigma^2}{\cos\beta}\left((x\cdot\theta_{\alpha})^2+(x\cdot\theta_{\beta})^2\cos^2\beta\right)\left(\theta\cdot\mathrm{d}x+(x\cdot\theta_{\alpha})\mathrm{d}\alpha+(x\cdot\theta_{\beta})\mathrm{d}\beta\right).
\end{equation}
So $\pi_R$ drops rank by 1 on $\Sigma=\{x\cdot\theta=0\}$ and its singularities are blowdown type as the kernel of $D\pi_R$ on $\Sigma$ is
\begin{equation}
\text{span}\left\{(x\cdot\theta_{\beta})\frac{\partial}{\partial\alpha}-(x\cdot\theta_{\alpha})\frac{\partial}{\partial\beta}\right\}\subset \ker\left\{\mathrm{d}\ \text{det} D\pi_R\right\}.
\end{equation}
Moreover, we have
\begin{equation}
\begin{split}
\pi_R(\Sigma)&=\{x\times\xi=0\}\\
&=\{(x,\xi) :p_i(x,\xi)=0, 1\leq i\leq 3\},
\end{split}
\end{equation}
where $p_1(x,\xi)=x_1\xi_2-x_2\xi_1$, $p_2(x,\xi)=x_1\xi_3-x_3\xi_1$ and $p_3(x,\xi)=x_2\xi_3-x_3\xi_2$. 

The Hamiltonian vector fields of the $p_i$ are given by
\begin{equation}
\begin{split}
H_{p_1}&=-x_2\partial_{x_1}+x_1\partial_{x_2}-\xi_2\partial_{\xi_1}+\xi_1\partial_{\xi_2},\\
H_{p_2}&=-x_3\partial_{x_1}+x_1\partial_{x_3}-\xi_3\partial_{\xi_1}+\xi_1\partial_{\xi_3},\\
H_{p_3}&=-x_3\partial_{x_2}+x_2\partial_{x_3}-\xi_3\partial_{\xi_2}+\xi_2\partial_{\xi_3}.
\end{split}
\end{equation}
Let $\rho=\sum_{i=1}^3\xi_i\partial_{\xi_i}$. Then, as $x=t \xi$ for some $t\in\mathbb{R}$, we can see that $\rho\notin\text{span}\{H_{p_1},H_{p_1},H_{p_3}\}$, so $\pi_R(\Sigma)$ is nonradial. 

To check that $\pi_R(\Sigma)$ is involutive, we first check that the Poisson brackets satisfy $\{p_i,p_j\}=0$, $i\neq j$:
\begin{equation}
\begin{split}
\{p_1,p_2\}&=H_{p_1}p_2=\xi_2x_3-x_2\xi_3=0\\
\{p_1,p_3\}&=H_{p_1}p_3=-\xi_1x_3+x_1\xi_3=0\\
\{p_2,p_3\}&=H_{p_2}p_3=\xi_1x_2-x_1\xi_2=0.
\end{split}
\end{equation}
Furthermore, if we work locally in a neighbourhood away from $x_1=0$, then $p_1,p_2=0\implies p_3=0$, so we need only consider the dependance of the differentials of $p_1$ and $p_2$. $\mathrm{d} p_1$ and $\mathrm{d} p_2$ are linearly independant if and only if the Hamiltonian vector fields of $p_1$ and $p_2$ are linearly independant. But $\text{span}\{H_{p_1},H_{p_2}\}$ has dimension 2, so $\pi_R(\Sigma)$ is involutive. So the conditions of Theorem \ref{blowdowntheorem} are satisfied except for the fact that $\pi_L$ is two-to-one away from $\Sigma$. However, we can remedy this by working locally in neighbourhoods of any given point $x_0$ within $B^3_{\alpha_1,\alpha_2}$ small enough so that $-x_0$ is not in the same neighbourhood. When we compose the operators restricted to neighbourhoods of $x_0$ and $-x_0$, we compose with the operator giving reflection in the origin so that even in that case we may apply Theorem \ref{blowdowntheorem}.

So, applying Theorem \ref{blowdowntheorem} we have $\mathcal{C}^*\mathcal{C}\in I^{2m+\frac{k-1}{2},-\frac{k-1}{2}}(\Delta \cup \widetilde{\Delta},\Lambda)$ where $k=1$ is the drop in rank of the left and right projections and $m=-1$ is the order of $\mathcal{C}$ as determined in Theorem \ref{FIO}.
\end{proof}
\end{theorem}
\noindent To complete this section we compute the flowout $\Lambda$ of $\pi_R(\Sigma)$.
\begin{corollary}
\label{corrflow}
Let $\pi_R$ be the right projection of $\mathcal{C}$ and let $\Sigma=\{x\cdot\theta =0\}$. Then the flowout of $\pi_R(\Sigma)$ is
\begin{equation}
\Lambda=\{\left(x,\xi;\mathcal{O}(x,\xi)\right) \ : \ x \in B^3_{\alpha_1,\alpha_2}, \quad \xi \in \mathbb{R}^3 \setminus 0, \quad x \times \xi = 0, \quad \mathcal{O} \in \Delta(SO_3 \times SO_3) \}.
\end{equation}
Here $\mathcal{O} (x,\xi) = (Ox, O\xi)$ where $O \in SO_3$ is any rotation.
\begin{proof}
Working locally away from $x_1=0$, we have
\begin{equation}
\begin{split}
\pi_R(\Sigma)=\{(x,\xi) :p_i(x,\xi)=0, 1\leq i\leq 2\},
\end{split}
\end{equation}
where $p_1(x,\xi)=x_1\xi_2-x_2\xi_1$ and $p_2(x,\xi)=x_1\xi_3-x_3\xi_1$. Letting $H_z=z_1H_{p_1}+z_2H_{p_2}$, by definition, the flowout of $\pi_R(\Sigma)$ is $\Lambda=\{(x,\xi;y,\eta)\in T^*X\times T^*X : (x,\xi)\in \pi_R(\Sigma), (y,\eta)=\exp(H_z)(x,\xi), z\in\mathbb{R}^2\}$.

We can write $H_z$ as
\begin{equation}
H_z=(x,\xi)\begin{pmatrix}
H^T & 0_{3\times 3}\\
0_{3\times 3} & H^T
\end{pmatrix}\begin{pmatrix}
\partial_{x}\\
\partial_{\xi}
\end{pmatrix},
\end{equation}
where $\partial_x=(\partial_{x_1},\partial_{x_2},\partial_{x_3})^T$, $\partial_{\xi}=(\partial_{\xi_1},\partial_{\xi_2},\partial_{\xi_3})^T$ and
\begin{equation}
H=\begin{pmatrix}
0 & -z_1 & -z_2\\
z_1 & 0 & 0\\
z_2 & 0 & 0
\end{pmatrix}.
\end{equation}
The flow of $H_z$ is thus given by the system of linear ODE's
\begin{equation}
\label{ODE}
\begin{pmatrix}
\dot{x}\\
\dot{\xi}
\end{pmatrix}=\begin{pmatrix}
H & 0_{3\times 3}\\
0_{3\times 3} & H
\end{pmatrix}\begin{pmatrix}
x\\
\xi
\end{pmatrix},
\end{equation}
with initial conditions $(x(0),\xi(0))=(x_0,\xi_0)$. Then the solution to (\ref{ODE}) at time $t'=1$ is
\begin{equation}
\begin{pmatrix}
{x}(1)\\
{\xi}(1)
\end{pmatrix}=\begin{pmatrix}
e^{H} & 0_{3\times 3}\\
0_{3\times 3} & e^{H}
\end{pmatrix}\begin{pmatrix}
x_0\\
\xi_0
\end{pmatrix},
\end{equation}
and hence the flow of $H_z$ can be computed as the exponential of the matrix $H$.

Now, if we parameterize $z_1=t\cos\omega$ and $z_2=t\sin\omega$ in terms of standard polar coordinates, then $H=tG=t(ba^T-ab^T)$, where $a=(1,0,0)^T$ and $b=(0,\cos\omega,\sin\omega)^T$. Let $P=-G^2=aa^T+bb^T$. Then $P^2=P$ ($P$ is idempotent) and $PG=GP=G$. From this it follows that
\begin{equation}
e^{H}=e^{tG}=I_{3\times 3}+G\sin t + G^2(1-\cos t).
\end{equation}
We can write
\begin{equation}
\hspace*{-0.4cm}
e^{H}=I_{3\times 3}+G\sin t + G^2(1-\cos t)=\begin{pmatrix}
\cos t & -\sin t\cos\omega & -\sin t\sin\omega\\
\sin t\cos\omega & \sin^2\omega+\cos t \cos^2\omega & -\cos t\cos\omega\sin\omega\\
\sin t\sin\omega & -\cos t\cos\omega\sin\omega & \cos^2\omega+\cos t \sin^2\omega 
\end{pmatrix},
\end{equation}
where $\omega\in [0,2\pi]$ and $t\in [0,\infty)$. It follows that
\begin{equation}
\label{rot}
e^{H}e^T_1=e^{H}(1,0,0)^T=(\cos t, \sin t\cos\omega, \sin t\sin\omega)^T.
\end{equation}
The right hand side of equation (\ref{rot}) is the standard parameterization of $S^2$ in terms of spherical coordinates, where $t$ is the polar angle from the $x$ axis pole and $\omega$ is the angle of rotation in the $yz$ plane (azimuth angle). So $e^H$ defines a full set of rotations on $S^2$ and hence, given a vector--conormal vector pair $(x,\xi)$, the Lagrangian $\Lambda$ includes the rotation of $x$ and $\xi$ over the whole sphere. This completes the proof.
\end{proof}
\end{corollary}

The above results tell us that the wavefront set of the kernel of the normal operator $\mathcal{C}^*\mathcal{C}$ is contained in $\Delta\cup \widetilde{\Delta} \cup\Lambda$, where $\Delta$ is the diagonal, $\widetilde{\Delta}$ the diagonal composed with reflection through the origin, and $\Lambda$ is the flowout from $\pi_R(\Sigma)$, which is a rotation by Corollary \ref{corrflow}. Also microlocally $\mathcal{C}^*\mathcal{C}\in I^{-2}(\Delta\backslash\Lambda)$ and $\mathcal{C}^*\mathcal{C}\in I^{-2}(\Lambda\backslash\Delta)$, which implies that the strength of the artefacts represented by $\Lambda$ are the same as the image intensity on $\Sigma=\{x\cdot \theta=0\}$. We will give examples of the artefacts implied by $\Lambda$ later in our simulations in section \ref{results}, but in the next section we shall show how to reduce the strength of this artefact microlocally.

\section{Reducing the strength of the image artifact}
\label{filtersection}
Here we derive a filter $Q$, which we show can be applied to reduce the strength of the image artefact $\Lambda$ for the cylinder tranform $\mathcal{C}$. We further show how $Q$ can be applied as a spherical convolution with a distribution $h$ on the sphere, which we will determine.

Using the ideas of \cite{felea}, our aim is to apply a filtering operator $Q:\mathcal{E}'((0,1) \times S^2)\to \mathcal{E}'((0,1) \times S^2)$, whose principal symbol vanishes to some order $s$ on $\pi_L(\Sigma)$, to $\mathcal{C}$ before applying the backprojection operator $\mathcal{C}^*$. From \cite{felea}, we have the following theorem.
\begin{theorem}
\label{filtertheorem}
Let $A\in I^m(C)$ be such that both the left projections of $A$, $\pi_L$ and $\pi_R$ are diffeomorphisms except on a set $\Sigma$ where they drop rank by $k$, and let $\pi_L(\Sigma)$ and $\pi_R(\Sigma)$ be involutive and nonradial. Let $Q$ be a pseudodifferential operator of order 0 whose principal symbol vanishes to order $s$ on $\pi_L(\Sigma)$. Then $A^*QA\in I^{2m+\frac{k-1}{2}-s,s-\frac{k-1}{2}}(\Delta,\Lambda)$, where $\Delta$ is the diagonal and $\Lambda$ is the flowout from $\pi_R(\Sigma)$.
\end{theorem}
Let $\Delta_{S^2}$ denote the Laplacian on $S^2$, and $I$ the identity operator. Then we will take
\begin{equation} \label{Qdef}
Q=-\Delta_{S^2}\left(I-\Delta_{S^2}\right)^{-1},
\end{equation}
whose symbol vanishes to order 2 on
\begin{equation}
\pi_L(\Sigma)=\{\hat{\alpha}=\hat{\beta}=0\}.
\end{equation}
There are two technical issues with the application of Theorem \ref{filtertheorem} in our case. One is the fact, which we already mentioned in the proof of Theorem \ref{maintheorem} that $\pi_L$ is two-to-one away from $\Sigma$. We can deal with this in the same way we dealt with in the proof of Theorem \ref{filtertheorem} by restricting $\mathcal{C}$ to small neighbourhoods of each point.

The other issue is that this operator $Q$, defined by \eqref{Qdef}, is not a pseudodifferential operator on $Y = (0,1) \times S^2$ since differentiation of its symbol in the dual angular variables does not increase the decay in the $\hat{s}$ direction. However, this objection can be overcome by noting that
\[
Q \underbrace{(1 - \Delta_{S^2})(1 - \Delta_{S^2} - \partial_s^2)^{-1}}_{\Psi_1}  = \underbrace{-\Delta_{S^2}(1 - \Delta_{S^2} - \partial_s^2)^{-1}}_{\Psi_2}.
\]
Both $\Psi_1$ and $\Psi_2$ are then pseudodifferential operators, and $\Psi_1$ is elliptic and of order zero. Thus, if $\Psi_1^{-1}$ is a pseudodifferential parametrix for $\Psi_1$, we have
\[
Q = \Psi_2 \Psi_1^{-1} + R
\] 
where $R$ is an operator with smooth kernel and $\Psi_2 \Psi_1^{-1}$ is a pseudodifferential operator satisfying the hypotheses of Theorem \ref{filtertheorem}. From this the results of Theorem \ref{filtertheorem} hold when $Q$ is given by \eqref{Qdef}.

We thus have, upon applying the filter $Q$ to $\mathcal{C}$ before applying $\mathcal{C}^*$ that $\mathcal{C}^*Q\mathcal{C}\in L^{-4,2}(\Delta, \Lambda)$. So $\mathcal{C}^*Q\mathcal{C}\in L^{-2}(\Delta\backslash\Lambda)$ and $\mathcal{C}^*Q\mathcal{C}\in L^{-4}(\Lambda\backslash\Delta)$ and the strength of the artefact is reduced and is now less than the strength of the image.

We now show how the filter $Q$ can be applied as a convolution with a distibution $h$ on the sphere. First we give some definitions and theorems on spherical harmonic expansions. For integers $l\geq 0, \ |m|\leq l$, we define the spherical harmonics $Y_l^m$ as
\begin{equation}
Y_l ^m(\alpha,\beta)=(-1)^m\sqrt{\frac{(2l+1)(l-m)!}{4\pi(l+m)!}}P_l ^m (\cos \beta) e^{\mathrm{i}m\alpha},
\end{equation}
where
\begin{equation}
P_l^m(x)=(-1)^m(1-x^2)^{m/2}\frac{\mathrm{d}^m}{\mathrm{d}x^m}P_l(x)
\end{equation}
and
\begin{equation}
P_l(x)=\frac{1}{2^l}\sum_{k=0}^{l}\binom{l}{k}^2(x-1)^{l-k}(x+1)^{k}
\end{equation}
are Legendre polynomials of degree $l$. The spherical harmonics $Y^m_l$ are the eigenfunctions of the Laplacian on $S^2$, with corresponding eigenvalues $c_l=-l(l+1)$. So $\Delta_{S^2}Y^m_l=c_lY^m_l$. From \cite{seeley} we have the following theorem.
\begin{theorem}
\label{series}
Let $F\in C^{\infty}(Z^3)$ and let
\begin{equation}
F_{lm}=\int_{S^2}F\bar{Y}_l ^m \mathrm{d}\Omega,
\end{equation}
where $\mathrm{d}\Omega$ is the surface measure on the sphere. Then the series
\begin{equation}
F_N=\sum_{0\leq l\leq N} \sum_{|m|\leq l}F_{lm}Y_{l}^m
\end{equation}
converges uniformly absolutely on compact subsets of $Z^3$ to $F$.
\end{theorem}
So after writing the cylinder tranform $\mathcal{C}$ in terms of its spherical harmonic expansion, we can apply the filter $Q$ as follows:
\begin{equation}
\label{applyQ}
\begin{split}
Q\mathcal{C}f(s,\theta)&=-\Delta_{S^2}\left(I-\Delta_{S^2}\right)^{-1}\sum_{l\in\mathbb{N}} \sum_{|m|\leq l}\mathcal{C}_{lm}(s)Y^m_l(\theta)\\
&=\sum_{l\in\mathbb{N}} \sum_{|m|\leq l}\frac{-c_l}{1-c_l}\mathcal{C}_{lm}(s)Y^m_l(\theta),
\end{split}
\end{equation}
where $\mathcal{C}_{lm}=\int_{S^2}\mathcal{C}\bar{Y}_l ^m \mathrm{d}\Omega$.

For a function $f$ on the sphere and $h$, a distribution on the sphere, we define the spherical convolution \cite{sphharm}
\begin{equation}
(f\ast_{S^2} h)(\theta)=\int_{g\in\text{SO}(3)}f(g\omega)h(g^{-1}\theta)\mathrm{d}g,
\end{equation}
where $\omega$ is the north pole, and from \cite{sphharm} we have the next theorem.
\begin{theorem}
\label{sphconv}
For functions $f,h\in L^2(S^2)$, the harmonic components of the convolution is a pointwise product of the harmonic components of the transforms:
\begin{equation}
(f\ast_{S^2} h)_{lm}=2\pi\sqrt{\frac{4\pi}{2l+1}}f_{lm}h_{l0}.
\end{equation}
\end{theorem}
\noindent For our case, this gives the following.
\begin{theorem}
\label{filter}
Let $Q=-\Delta_{S^2}\left(I-\Delta_{S^2}\right)^{-1}$ and let $F\in C^{\infty}_0(Z^3)$. Then
\begin{equation}
QF=\left(F\ast_{S^2} h\right),
\end{equation}
where $h$ is defined by
\begin{equation}
h(\theta)=\sum_{l\in\mathbb{N}}h_l\sum_{|m|\leq l}Y^m_l(\theta),
\end{equation}
where
\begin{equation}
h_l=\frac{l(l+1)a_l}{l(l+1)+1}
\end{equation}
and $a_l=\frac{1}{2\pi}\cdot\sqrt{\frac{2l+1}{4\pi}}$.
\begin{proof}
From equation (\ref{applyQ}) we have
\begin{equation}
\begin{split}
QF(s,\theta)&=\sum_{l\in\mathbb{N}} \sum_{|m|\leq l}\frac{l(l+1)}{l(l+1)+1}F_{lm}(s) Y^m_l(\theta)\\
&=\sum_{l\in\mathbb{N}} \sum_{|m|\leq l}2\pi\sqrt{\frac{4\pi}{2l+1}}h_lF_{lm}(s) Y^m_l(\theta).
\end{split}
\end{equation}
Defining $h_{lm}=h_l$ for all $l\in \mathbb{N}$, $|m|\leq l$, we have by Theorem \ref{sphconv}
\begin{equation}
\begin{split}
QF(s,\theta)&=\int_{g\in\text{SO}(3)}F(s,g\omega)h(g^{-1}\theta)\mathrm{d}g\\
&=\left(F\ast_{S^2} h\right)(s,\theta),
\end{split}
\end{equation}
where $h(\theta)=\sum_{l\in\mathbb{N}} \sum_{|m|\leq l}h_{lm}Y^m_l(\theta)=\sum_{l\in\mathbb{N}}h_l\sum_{|m|\leq l}Y^m_l(\theta)$, which completes the proof.
\end{proof}
\end{theorem}
\section{Simulations}
\label{results}
Given the equivalence of the spindle transform $\mathcal{S}$ and the cylinder transform $\mathcal{C}$, and given also that the diffeomorphism defining their equivalence $v(x)=\left(\sqrt{1+\frac{1}{|x|^2}}-\frac{1}{|x|}\right)\cdot \frac{x}{|x|}$ depends only on $|x|$, the artefacts described by the Lagrangian $\Lambda$ apply also to the normal operator of the spindle transform $\mathcal{S}^*\mathcal{S}$. Here we simulate the image artefacts produced by $\Lambda$ in image reconstructions from spindle transform data and show how the filter we derived in section \ref{filtersection} can be used to reduce these artefacts. We also provide simulated reconstructions of densities which we should find difficult to reconstruct from a microlocal perspective  (i.e. densities whose wavefront set is in directions normal to the surface of a sphere centred at the origin), and investigate the effects of applying the filter $Q$ as a pre--conditioner, prior to implementing some discrete solver, in our reconstruction. 

To conduct our simulations we consider the discrete form of the spindle transform as in \cite{WL}, and solve the linear system of equations
\begin{equation}
Ax=b,
\end{equation}
where $A$ is the discrete operator of the spindle transform, $x$ is the vector of pixel values and $b=Ax$ is the vector of spindle transform values (simulated as an inverse crime). To apply the filter $Q$ derived in section \ref{filtersection}, we decompose $b$ into its first $L$ spherical harmonic components and then multiply each component by the filter components $h_l/a_l=\frac{l(l+1)}{l(l+1)-1}$, $0\leq l\leq L$ before recomposing the series.

\begin{figure}[h]
\hspace*{-2cm}
\begin{subfigure}{0.5\textwidth}
\includegraphics[width=1.2\linewidth, height=10.5cm]{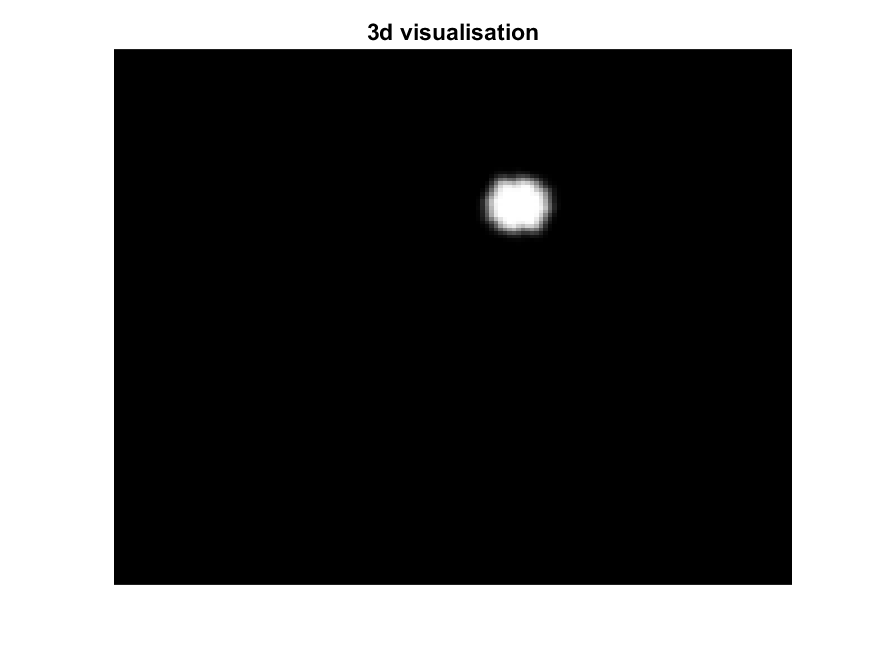} 
\end{subfigure}\hspace{15mm}
\begin{subfigure}{0.5\textwidth}
\includegraphics[width=1.2\linewidth, height=10.5cm]{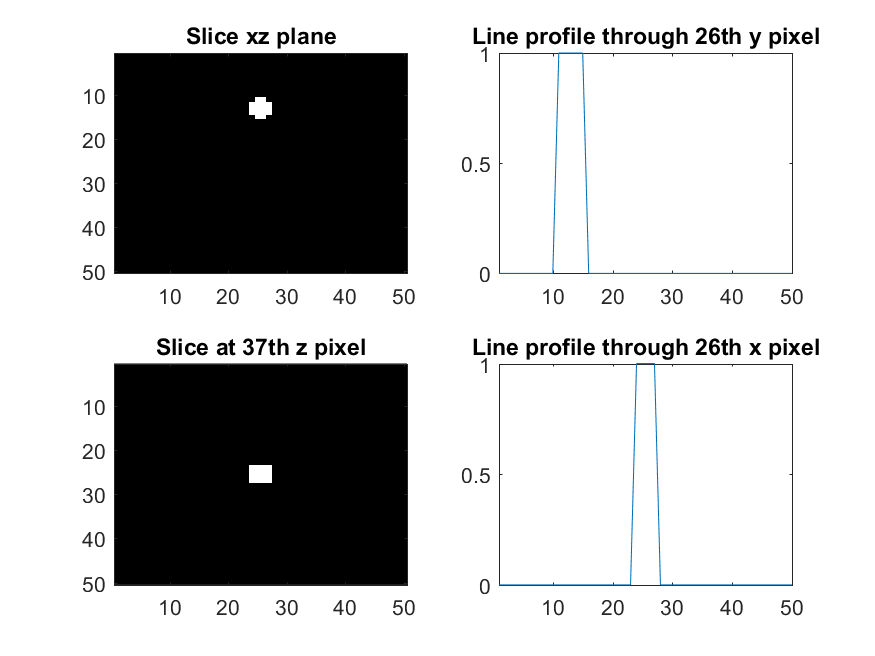}
\end{subfigure}
\caption{Small bead.}
\label{F1}
\hspace*{-2cm}
\begin{subfigure}{0.5\textwidth}
\includegraphics[width=1.2\linewidth, height=10.5cm]{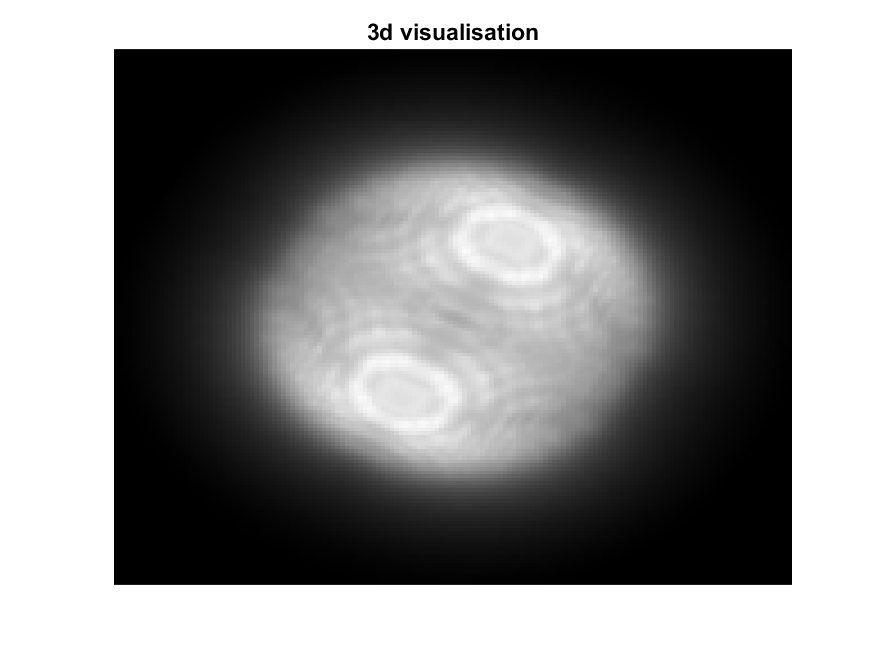} 
\end{subfigure}\hspace{15mm}
\begin{subfigure}{0.5\textwidth}
\includegraphics[width=1.2\linewidth, height=10.5cm]{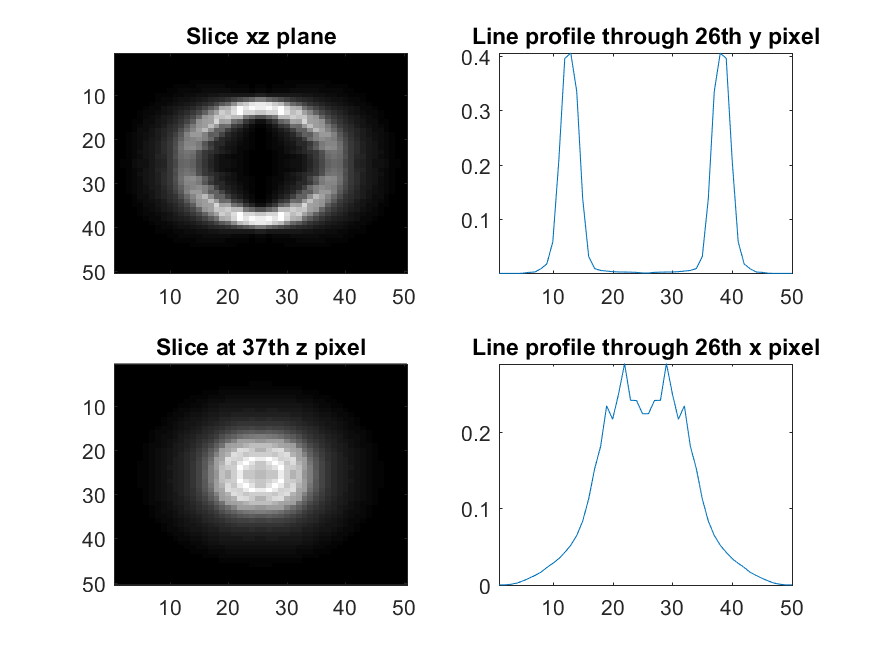}
\end{subfigure}
\caption{Bead reconstruction by backprojection.}
\label{F2}
\end{figure}

\begin{figure}[h]
\hspace*{-2cm}
\begin{subfigure}{0.5\textwidth}
\includegraphics[width=1.2\linewidth, height=10.5cm]{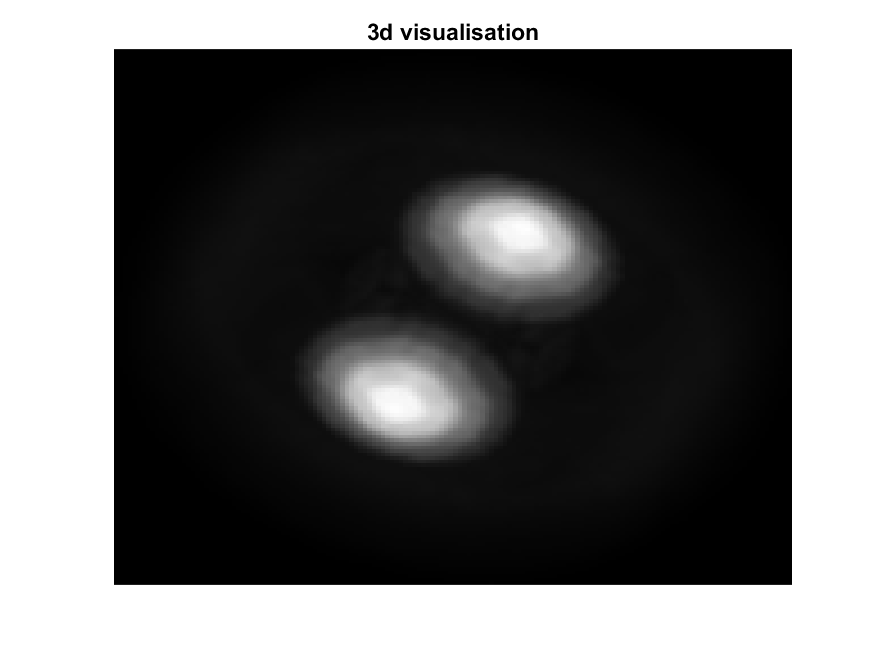} 
\end{subfigure}\hspace{15mm}
\begin{subfigure}{0.5\textwidth}
\includegraphics[width=1.2\linewidth, height=10.5cm]{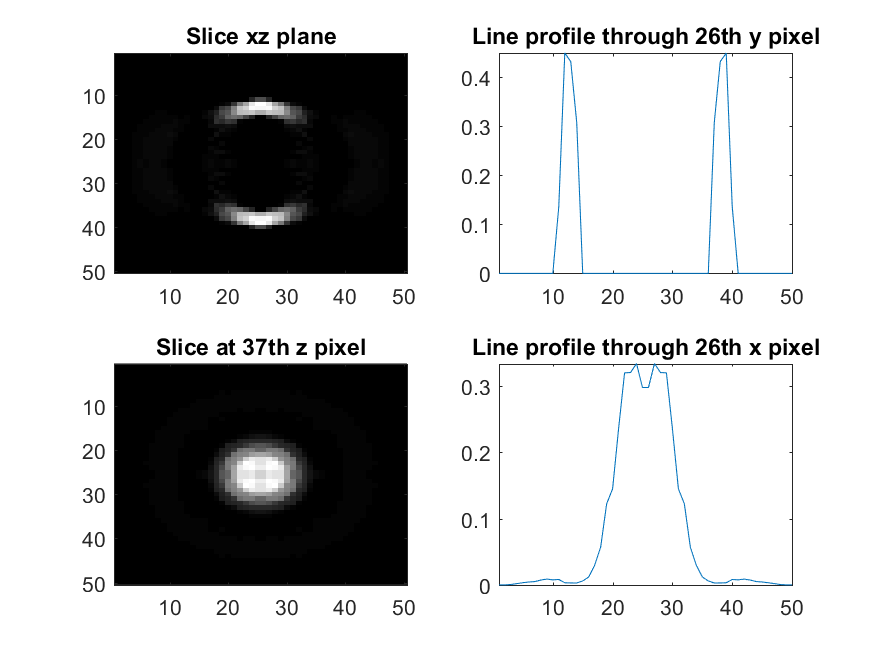}
\end{subfigure}
\caption{Bead reconstruction by filtered backprojection, with $L=25$ components.}
\label{F3}
\hspace*{-2cm}
\begin{subfigure}{0.5\textwidth}
\includegraphics[width=1.2\linewidth, height=10.5cm]{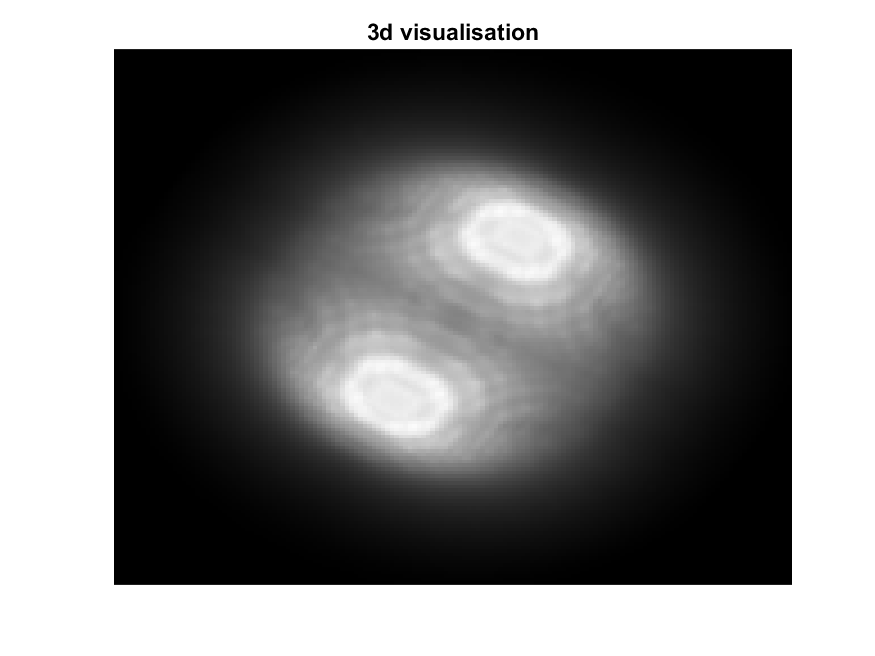} 
\end{subfigure}\hspace{15mm}
\begin{subfigure}{0.5\textwidth}
\includegraphics[width=1.2\linewidth, height=10.5cm]{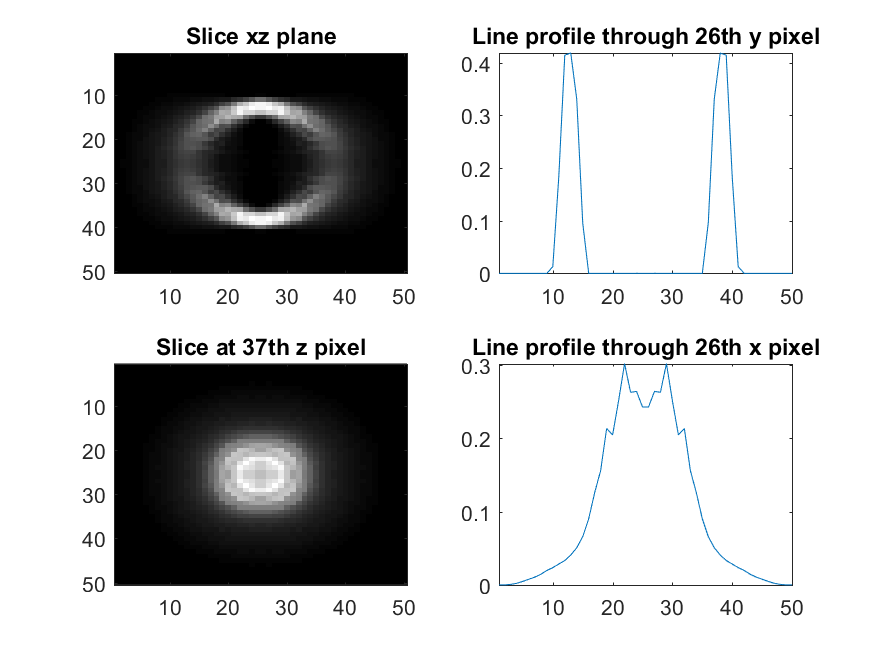}
\end{subfigure}
\caption{Bead reconstruction by backprojection, truncating the data to $L=25$ components.}
\label{F4}
\end{figure}

Consider the small bead of constant density pictured in figure \ref{F1}. In figure \ref{F2} we present a reconstruction of the small bead by unfiltered backprojection (represented as an MIP image to highlight the bead). Here we see artefacts described by the Lagrangian $\Lambda$ as, in the reconstruction, the bead is smeared out over the sphere. If we apply the filter $Q$ to the spherical components of the data and sum over the first $L=25$ components before backprojecting, then we see a significant reduction in the strength of the artefact, the density is more concentrated around the small bead and the image is sharper. See figure \ref{F3}. If we simply truncate the harmonic series of our data before backprojecting without a filter, this has a regularising effect and the level of blurring around the sphere is reduced. However we still see the artefacts due to $\Lambda$. See figure \ref{F4}. The line profiles in figures \ref{F1}--\ref{F4} have been normalised.
We note that in the reconstructions presented, the object is reflected through the origin in the reconstruction. This is as predicted by the Lagrangian $\widetilde{\Delta}$. The reflection artefact seems intuitive given the symmetries involved in our geometry. It was shown in \cite{WL} that the null space of $\mathcal{S}$ consists of odd functions (i.e. functions whose even harmonic components are zero), and so what we see in the reconstruction is the projection of the density onto its even components. We also see this effect in the reconstructions presented in \cite{WL}.

\begin{figure}[h]
\hspace*{-2cm}
\begin{subfigure}{0.5\textwidth}
\includegraphics[width=1.2\linewidth, height=10.5cm]{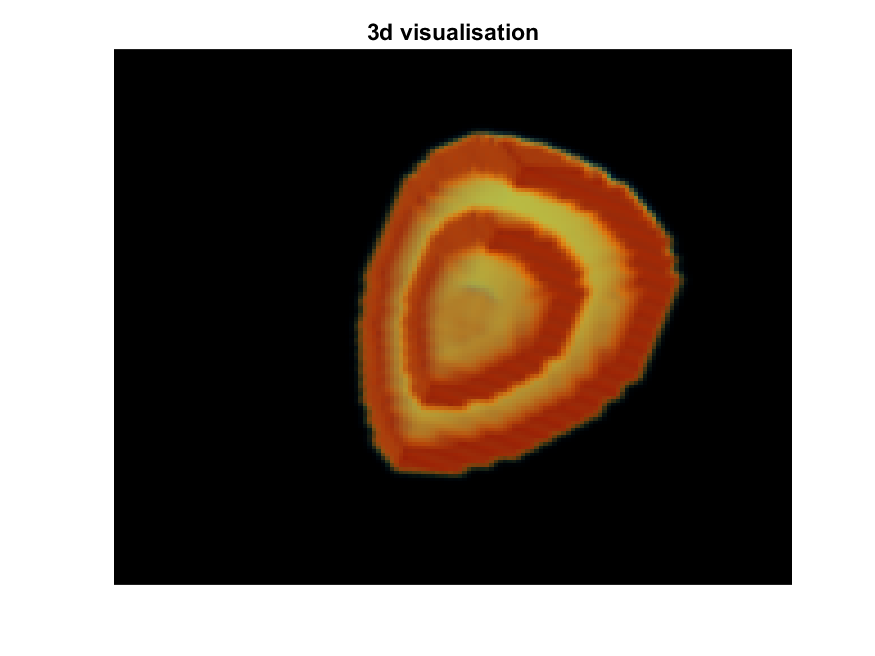} 
\end{subfigure}\hspace{15mm}
\begin{subfigure}{0.5\textwidth}
\includegraphics[width=1.2\linewidth, height=10.5cm]{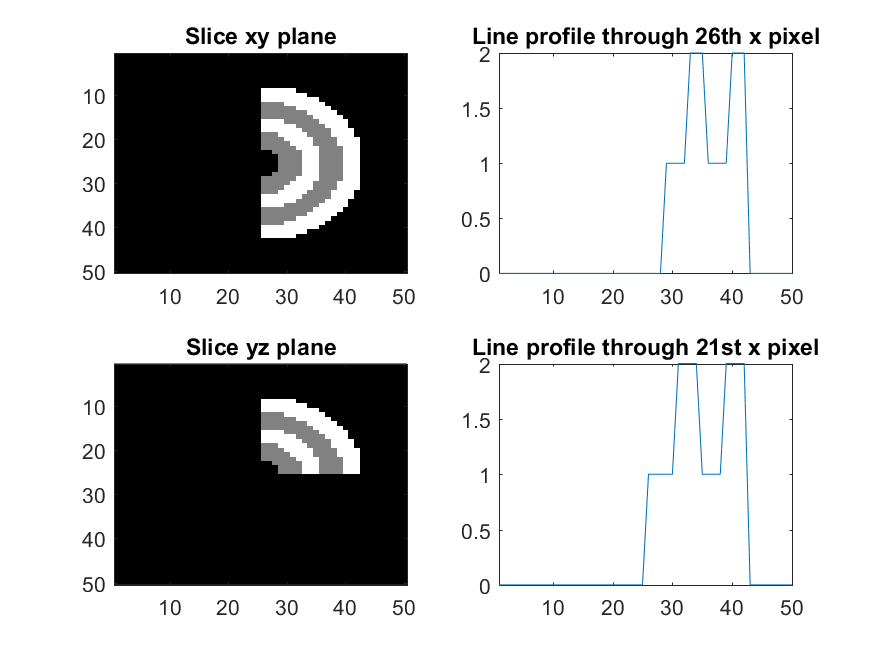}
\end{subfigure}
\caption{Layered spherical shell segment phantom, centred at the origin.}
\label{F5}
\hspace*{-2cm}
\begin{subfigure}{0.5\textwidth}
\includegraphics[width=1.2\linewidth, height=10.5cm]{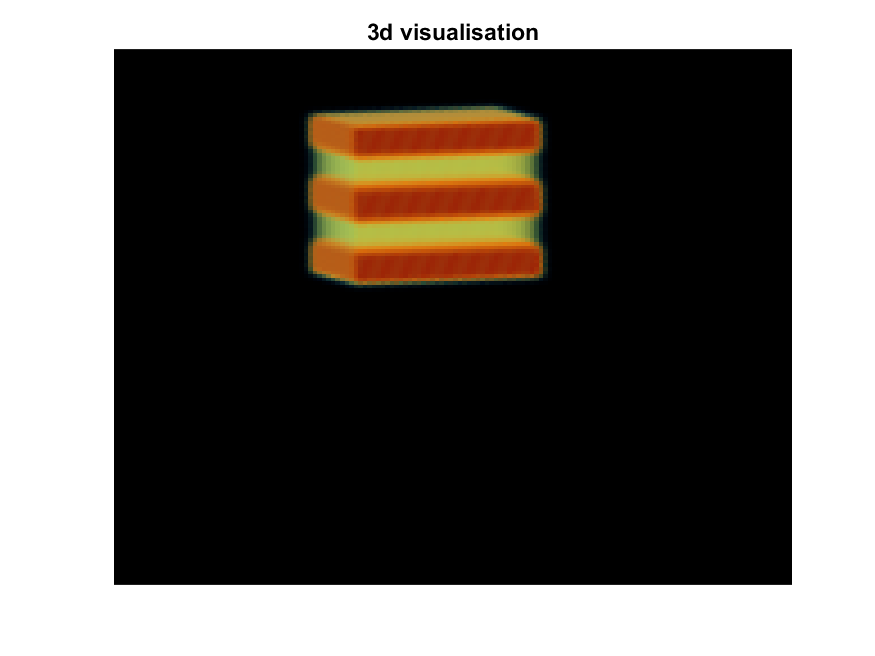} 
\end{subfigure}\hspace{15mm}
\begin{subfigure}{0.5\textwidth}
\includegraphics[width=1.2\linewidth, height=10.5cm]{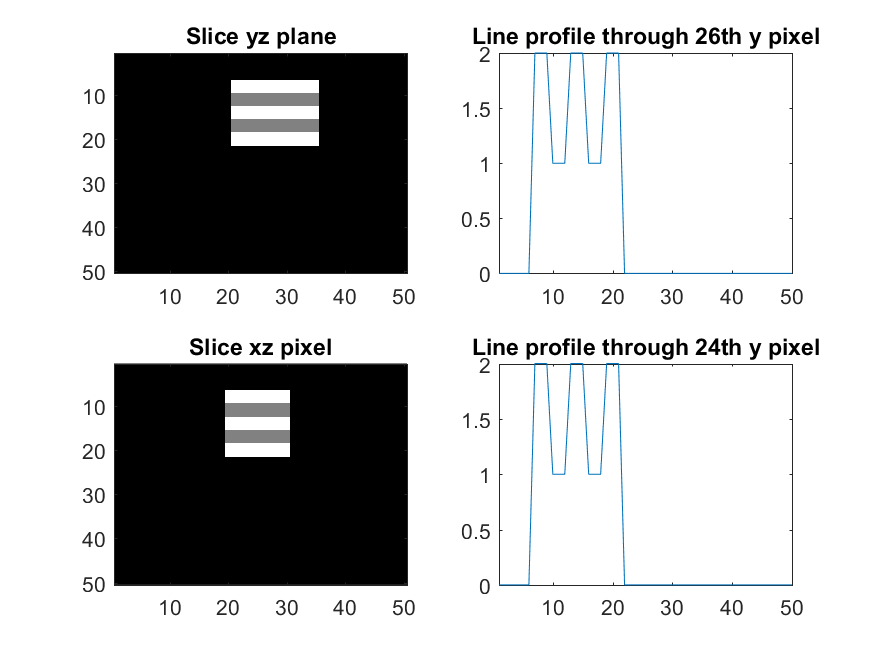}
\end{subfigure}
\caption{Layered plane phantom.}
\label{F7}
\end{figure}

\begin{figure}[h]
\hspace*{-2cm}
\begin{subfigure}{0.5\textwidth}
\includegraphics[width=1.2\linewidth, height=10.5cm]{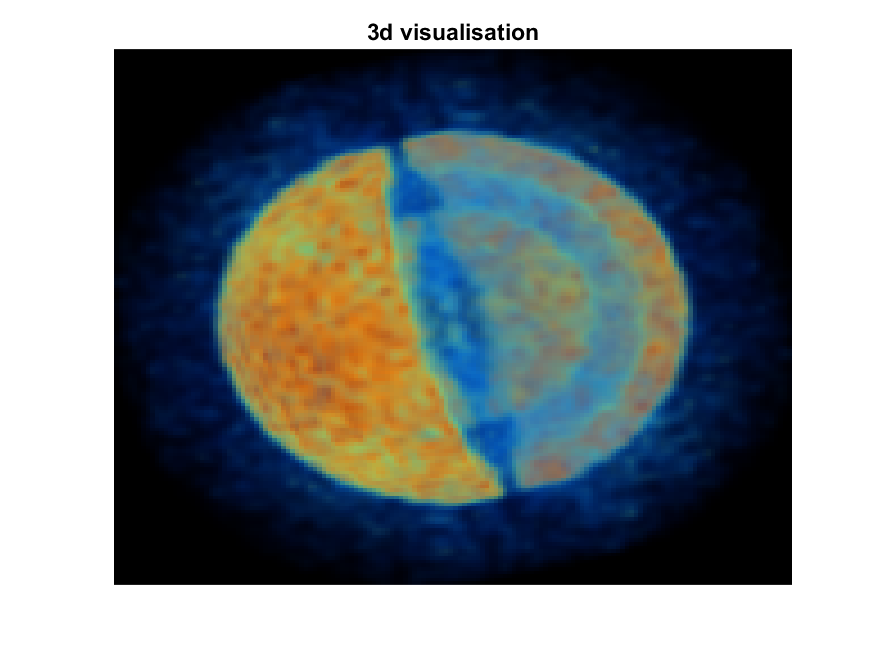} 
\end{subfigure}\hspace{15mm}
\begin{subfigure}{0.5\textwidth}
\includegraphics[width=1.2\linewidth, height=10.5cm]{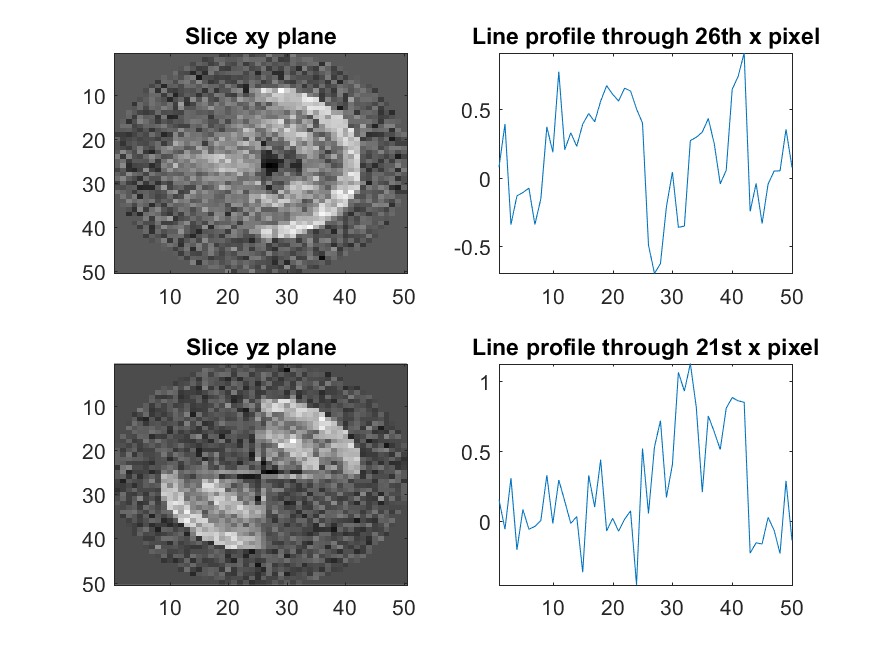}
\end{subfigure}
\caption{Layered spherical shell segment CGLS reconstruction.}
\label{F6}
\hspace*{-2cm}
\begin{subfigure}{0.5\textwidth}
\includegraphics[width=1.2\linewidth, height=10.5cm]{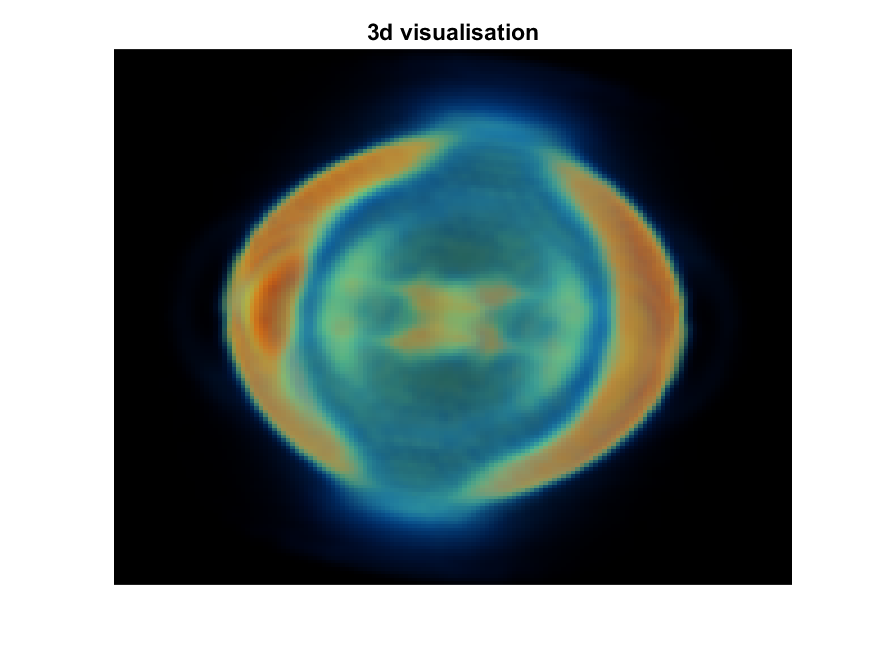} 
\end{subfigure}\hspace{15mm}
\begin{subfigure}{0.5\textwidth}
\includegraphics[width=1.2\linewidth, height=10.5cm]{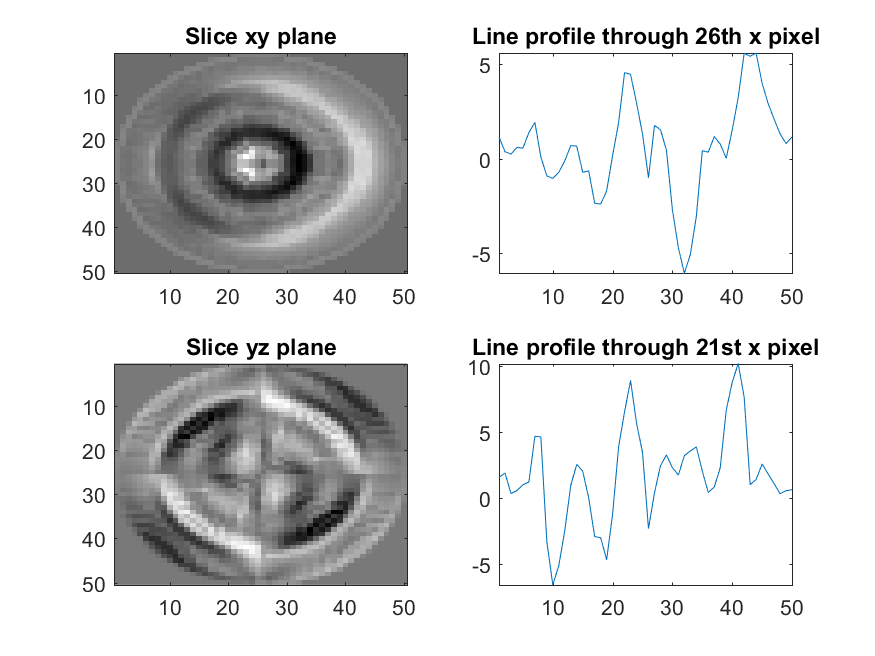}
\end{subfigure}
\caption{Layered spherical shell segment CGLS reconstruction, with $Q^{\frac{1}{2}}$ used as a pre--conditioner and no added Tikhonov regularisation.}
\label{F10}
\end{figure}

\begin{figure}[h]
\hspace*{-2cm}
\begin{subfigure}{0.5\textwidth}
\includegraphics[width=1.2\linewidth, height=10.5cm]{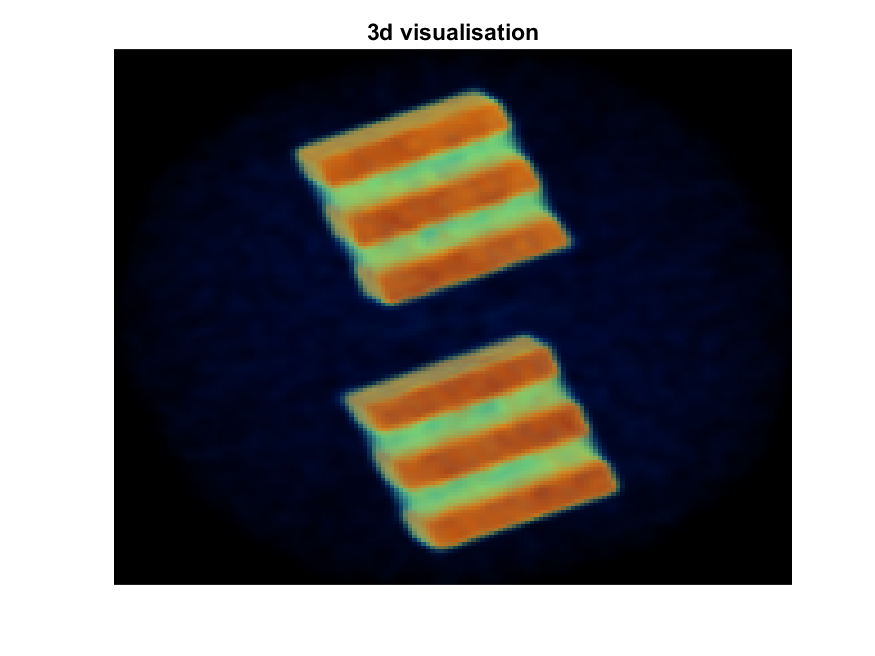} 
\end{subfigure}\hspace{15mm}
\begin{subfigure}{0.5\textwidth}
\includegraphics[width=1.2\linewidth, height=10.5cm]{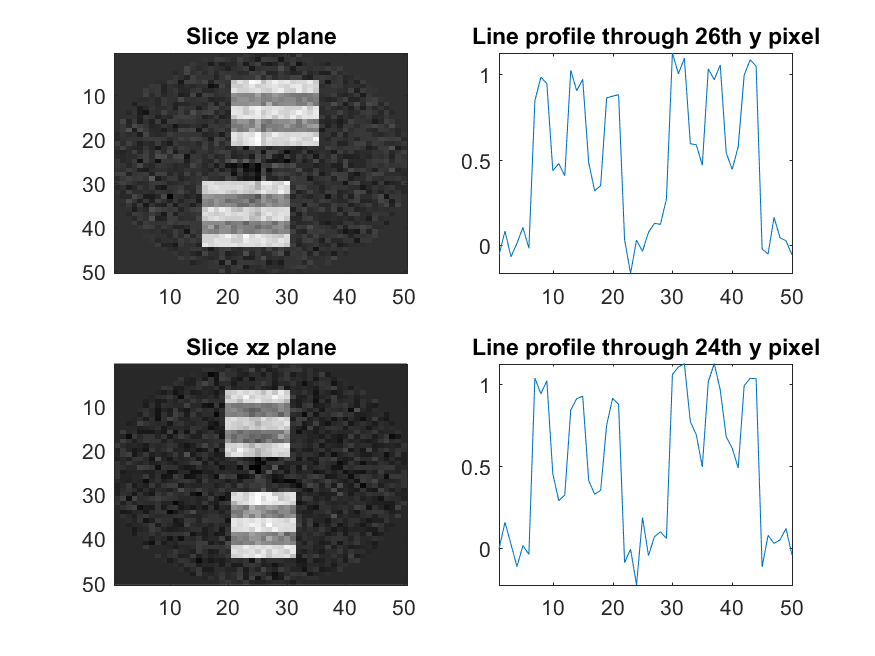}
\end{subfigure}
\caption{Layered plane reconstruction by CGLS.}
\label{F8}
\hspace*{-2cm}
\begin{subfigure}{0.5\textwidth}
\includegraphics[width=1.2\linewidth, height=10.5cm]{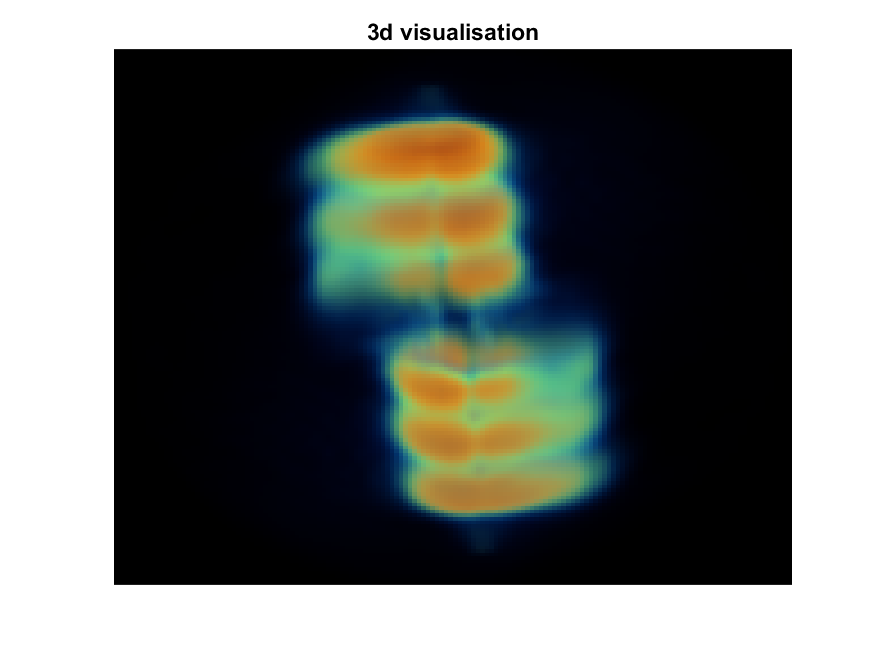} 
\end{subfigure}\hspace{15mm}
\begin{subfigure}{0.5\textwidth}
\includegraphics[width=1.2\linewidth, height=10.5cm]{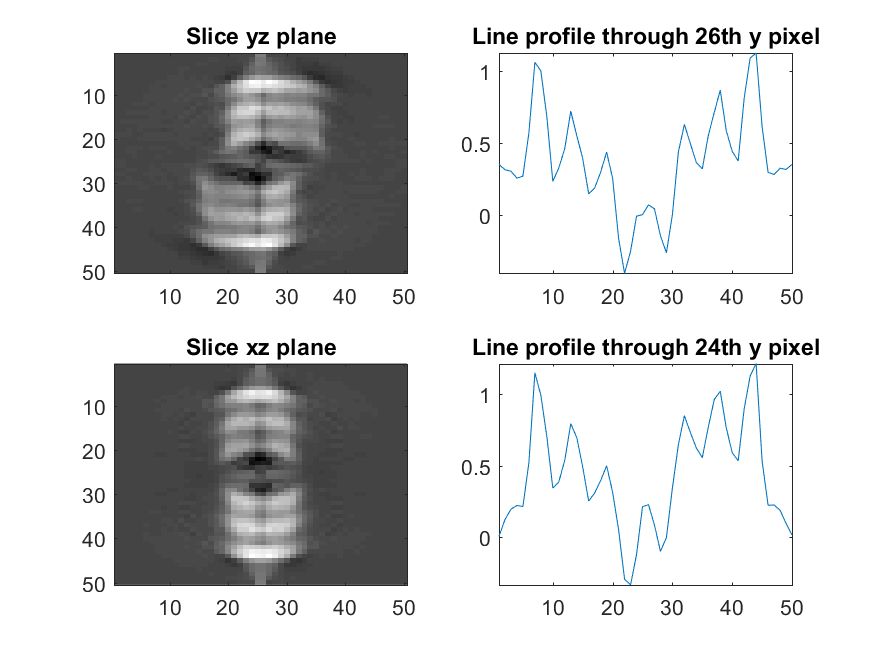}
\end{subfigure}
\caption{Layered plane reconstruction by Landweber iteration.}
\label{F11}
\end{figure}

Now let us consider the layered spherical shell segment phantom (the layers have values oscillating between 1 and 2) centred at the origin, shown in figure \ref{F5}. We reconstruct the phantom by applying CGLS implicitly to the normal equations (i.e. we avoid a direct application of $A^TA$) with $1\%$ added Gaussian noise and regularise our solution using Tikhonov regularisation. See figure \ref{F6}. Here the image quality is not clear and the layers seem to blur into one, and the jump discontinuities in the image are not reconstructed adequately. However if we arrange the layers as sections of planes and perform the same reconstruction (see figures \ref{F7} and \ref{F8}), then the image quality is significantly improved, and the jump discontinuities between the oscillating layers are clear. This is as expected, as the wavefront set of the  spherical density is contained in $\pi_R(\Lambda)$, so we see artefacts in the reconstruction. When the layers are arranged as planes this is not the case and we see an improvement in the reconstruction. In figure \ref{F10} we have investigated the effects of applying the filter $Q$ as a pre--conditioner prior to a CGLS implementation. To obtain the reconstruction, we solved the system of equations $Q^{\frac{1}{2}}Ax=Q^{\frac{1}{2}}b$ using CGLS with $1\%$ added Gaussian noise. The filter has the effect of smoothing the radial singularities in the reconstruction. Here we see that the outer shell is better distinguished than before but the inner shells fail to reconstruct and overall the image quality is not good. 

\begin{figure}[h]
\hspace*{-2cm}
\begin{subfigure}{0.5\textwidth}
\includegraphics[width=1.2\linewidth, height=10.5cm]{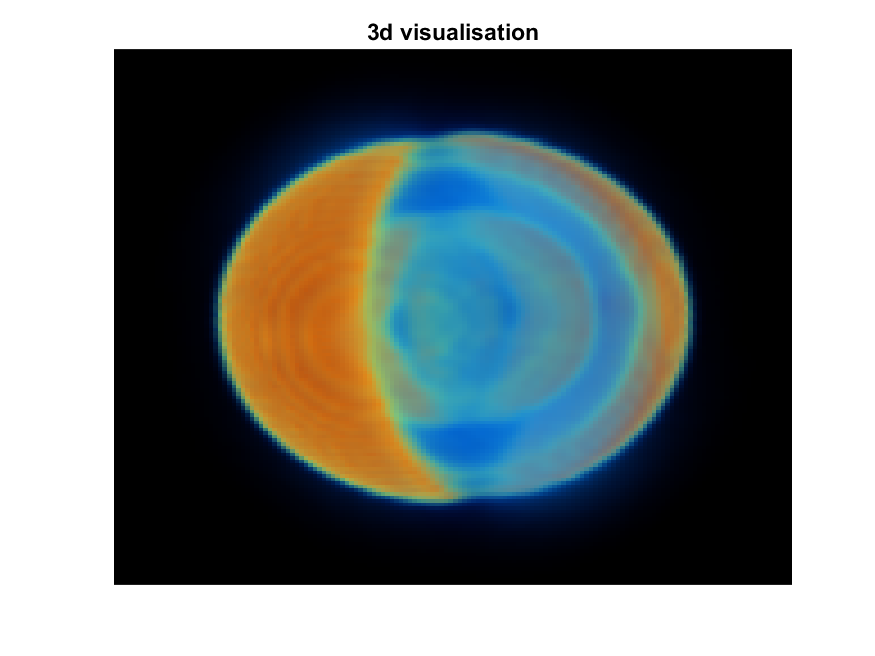} 
\end{subfigure}\hspace{15mm}
\begin{subfigure}{0.5\textwidth}
\includegraphics[width=1.2\linewidth, height=10.5cm]{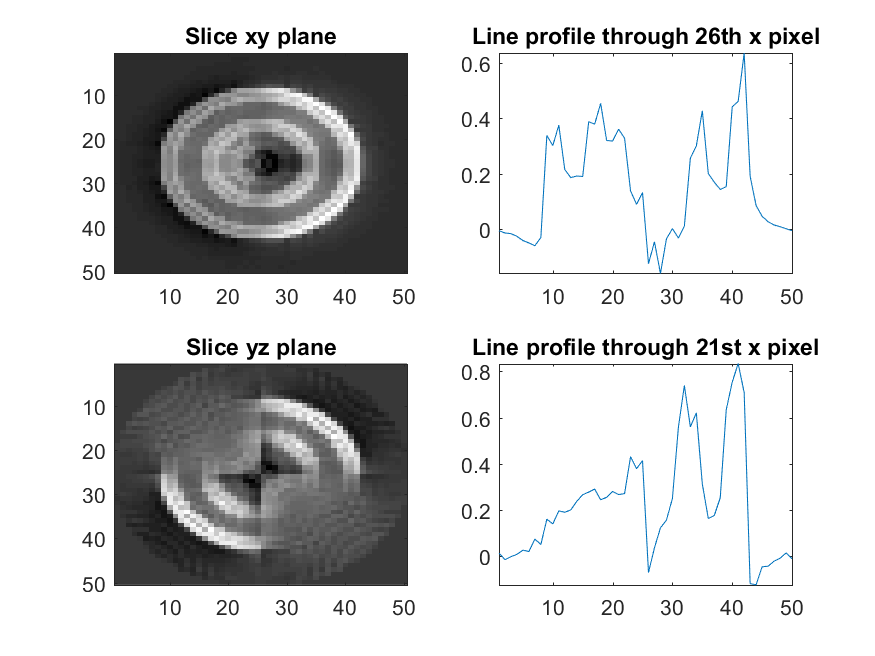}
\end{subfigure}
\caption{Spherical shell reconstruction by Landweber iteration.}
\label{F12}
\hspace*{-2cm}
\begin{subfigure}{0.5\textwidth}
\includegraphics[width=1.2\linewidth, height=10.5cm]{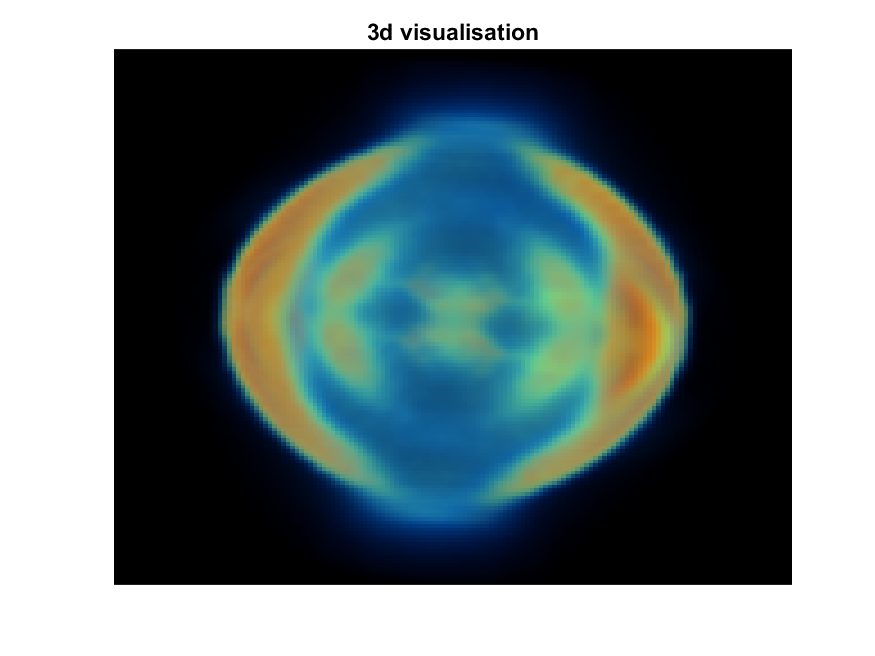} 
\end{subfigure}\hspace{15mm}
\begin{subfigure}{0.5\textwidth}
\includegraphics[width=1.2\linewidth, height=10.5cm]{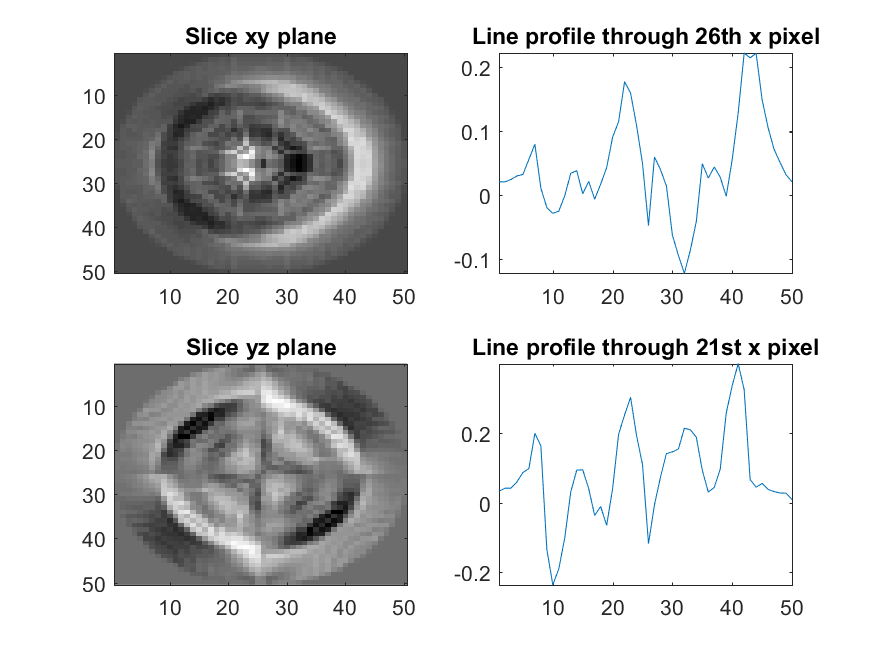}
\end{subfigure}
\caption{Spherical shell reconstruction by Landweber iteration, with $Q^{\frac{1}{2}}$ used as a pre--conditioner.}
\label{F13}
\end{figure}

In figures \ref{F11}, \ref{F12} and \ref{F13} we have presented reconstructions of the layered spherical shell and layered plane phantoms by Landweber iteration, with $1\%$ Gaussian noise. Here the jump discontinuities in the spherical shell reconstruction are clearer. However as the Landweber method applies the normal operator ($A^TA$) at each iteration, we see the artefacts predicted by $\Lambda$ in the reconstruction and the spherical segment is blurred out over spheres centred at the origin. The artefacts are less prevalent in the plane phantom reconstruction. Although we do see some blurring at the plane edges. In figure \ref{F13} we have applied $Q^{\frac{1}{2}}$ as a pre--conditioner to a Landweber iteration. Here it is not clear that we see a reduction in the spherical artefact and there is a loss in clarity due to the level of smoothing.
\clearpage
\clearpage
\section{Conclusions and further work}
We have presented a microlocal analysis of the spindle transform introduced in \cite{WL}. An equivalence to a cylinder transform $\mathcal{C}$ was proven and the microlocal properties of $\mathcal{C}$ were studied. We showed that $\mathcal{C}$ was an FIO whose normal operator belonged to a class of distributions $I^{p,l}(\Delta,\Lambda)$, where $\Lambda$ is the flowout from the right projection of $\mathcal{C}$, which we calculated explicitly. In section \ref{filtersection}, we showed how to reduce the size of the rotation artefact associated to $\Lambda$ microlocally, through the application of an operator $Q$, and showed that $Q$ could be applied as a spherical convolution with a distribution $h$ on the sphere, or using spherical harmonics. We provided simulated reconstructions to show the artefacts produced by $\Lambda$, and showed how applying $Q$ reduced the artefacts in the reconstruction. Reconstructions of densities of oscillating layers were provided using CGLS and a Landweber iteration. We also gave reconstructions of a spherical layered shell centred at the origin, using $Q^{\frac{1}{2}}$ as a pre--conditioner, prior to a CGLS and Landweber implementation and compared our results. 

In future work we aim to derive an inversion formula of either a filtered backprojection or backprojection filter type. That is, we aim to determine whether there exists an operator $A$ such that either $\mathcal{S}^*A$ or $A\mathcal{S}^*$ is a left inverse for $\mathcal{S}$. After which we could see how the filter $Q$ derived here may be involved in the inversion process. We also aim to assess if Sobolev space estimates can be derived for the spindle transform to gain a further understanding of its stability.

\section{Acknowledgements}

The authors would like to thank Bill Lionheart for suggesting the project, and for useful discussions on the topic.

\end{document}